\documentclass[10pt,a4paper]{article}

\usepackage[latin1]{inputenc}
\usepackage{amsmath}
\usepackage{amsfonts}
\usepackage{amssymb}
\usepackage{multirow}
\usepackage{colortbl}
\usepackage{rotating}
\usepackage{caption}
\usepackage{verbatim}
\usepackage{algorithmic}
\usepackage{algorithm}
\usepackage{paralist}

\usepackage{times}

\DeclareMathOperator{\lcm}{lcm}
\DeclareMathOperator{\smallest}{smallest}

\def \refsec {\S}

\newtheorem{lemma}{Lemma}

\begin{document}

\author{
Arthur R. Townsend
\\
 \texttt{mail@arthurtownsend.com}
}

\title{The Search for Maximal Values of 
$\min(A,B,C)/ \gcd(A,B,C)$ for $A^x + B^y = C^z$}

\date{\today}

\maketitle

\let\oldthebibliography=\thebibliography
  \let\endoldthebibliography=\endthebibliography
  \renewenvironment{thebibliography}[1]{%
    \begin{oldthebibliography}{#1}%
      \setlength{\parskip}{0ex}%
      \setlength{\itemsep}{0ex}%
  }%
  {%
    \end{oldthebibliography}%
  }

\begin{abstract}

This paper answers a question asked by Ed Pegg~Jr. in 2001: ``What is the maximal value of
$\min(A,B,C) / \gcd(A,B,C)$ for $A^x + B^y = C^z$ with $A,B,C \geq 1;$ $x,y,z \geq 3$?''
Equations of this form are analyzed, showing how they map to exponential Diophantine 
equations with coprime bases. 
A search
algorithm is provided to find the largest min/gcd value within a given equation range.  
The algorithm precalculates a multi-gigabyte lookup table of power residue information that is used to eliminate
over 99\% of inputs with a single array lookup and without any further calculations.  On inputs that pass this test,
the algorithm then performs further power residue tests, 
avoiding modular powering by using lookups into precalculated
tables, and avoiding division by using multiplicative inverses.  This algorithm is used to show
the largest min/gcd value for all equations with $C^z \leq 2^{100}$.
\end{abstract}

\noindent
The Tijdeman-Zagier Conjecture (T-Z Conjecture)
\cite{bib-bruin-sum-two-cubes,bib-crandall,bib-elkies,bib-waldschmidt},
also known as the Beal Prize Problem\cite{bib-mauldin},
is the claim that the following equation has no solutions: 

Tijdeman-Zagier Conjecture Equation / Beal Prize Problem Equation

$A^x + B^y = C^z$
\qquad
$x,y,z \in \mathbb{Z} _{\geq 3}$
\qquad
$A,B,C \in \mathbb{Z} ^+$
\qquad
$\gcd(A,B,C) = 1$
\\
Regarding this conjecture,
in 2001,
Ed Pegg~Jr. asked the following question
\cite{bib-pegg}:
Without the restriction that
$A,B,C$ be coprime,
what is the maximal value of $\min(A,B,C)$ / $\gcd(A,B,C)$?
We define the \emph{Pegg Value} to be $\min(A,B,C) / \gcd(A,B,C)$ of a T-Z Conjecture Equation or
Resultant Pegg Equation.

Resultant Pegg Equation

$A^x + B^y = C^z$
\qquad
$x,y,z \in \mathbb{Z} _{\geq 3}$
\qquad
$A,B,C \in \mathbb{Z} ^+$
\qquad
$\gcd(A,B,C) > 1$
\\
Notice that in the T-Z Conjecture, the $\gcd(A,B,C) = 1$ condition avoids an infinite 
number of trivial and uninteresting solutions.  For example, adding any two 
integers $A$ and $B$ gives an equation,
$A + B = C$
which multiplied by
$A^{20} B^{15} C^{24}$
yields
$(A^7 B^5 C^8)^3 + (A^5 B^4 C^6)^4 = (A^4 B^3 C^5)^5$.
The Pegg Value for this construction is 1, even if the original terms are all coprime.  In contrast 
to this construction, any T-Z Conjecture counterexample will have a Pegg Value equal to 
$\min(A,B,C).$ And $\min(A,B,C)$ will
be $ > 1$, as Preda Mih\u{a}ilescu showed\cite{bib-mihailescu} that the only perfect powers that differ by 1
are $2^3$ and $3^2$.   So in some sense, equations with a Pegg Value $> 1$ are ``closer'' to 
being T-Z Conjecture counterexamples.
Let us examine precisely what types of equations have a Pegg Value $> 1$ that are not T-Z 
Conjecture counterexamples.

\section{Original Equations and Conversion to Resultant Equations}
\label{section1}

Starting with a Resultant Pegg Equation 
with a Pegg Value $> 1$,
then expanding and dividing out the common factors 
(so the terms to be added are coprime), 
one is left with:
$d a^x + e b^y = f c^z;$
$x,y,z \geq 3;$
$a,b,c,d,e,f \geq 1;$
$\gcd(d a, e b, f c) = 1;$
and there exists a number $N$ such that $N d$ is an $x$-th power, $N e$ is a $y$-th 
power, and $N f$ is a $z$-th power.
If the smallest $N = 1$,
then the Resultant Pegg Equation is simply a T-Z Conjecture Counterexample that has
been multiplied by a number that is simultaneously an $x$-th power, $y$-th power, and $z$-th power.
Otherwise
at least one of $d,e,f > 1,$
and we call the equation an \emph{Original Pegg Equation}.
Multiplying this original equation
$d a^x + e b^y = f c^z$
by $N$ yields its resultant equation:
$A^x + B^y = C^z$
with
$A = (N d)^{1/x}$ $a;$
$B = (N e)^{1/y}$ $b;$
$C = (N f)^{1/z}$ $c\enspace.$

For ease of reference throughout this paper, $a,b,c$ will refer to the bases to the 
$x,y,z$ exponents, respectively, of the original equation.
And $A,B,C$ will refer to the bases to the 
$x,y,z$ exponents, respectively, of the resultant equation.
Similarly, $d,e,f$ will refer to the original coefficients associated with the $a,b,c$ bases,
respectively.  And $D,E,F$ will refer to the resultant coefficients associated with the 
resultant $A,B,C$ bases, respectively. $D,E,F$ are the portions of the resultant 
bases that result from ``spreading'' the original coefficients across the equation to make 
each term a perfect power:
$D = (N d)^{1/x} = A / a;$
$E = (N e)^{1/y} = B / b;$
$F = (N f)^{1/z} = C / c$.
So the original equation
$d a^x + e b^y = f c^z$
after multiplying by N, converts to
$(D a)^x + (E b)^y = (F c)^z$
or equivalently
$A^x + B^y = C^z\enspace.$
If $a,b,c,d,e,f$ are all pairwise coprime we see that the Pegg Value cannot be less than $\min(a,b,c),$
because one of $a,b,c$ will be a factor of $\min(A,B,C)$ and $\gcd(A,B,C)$ will not contain any 
factors of $a,b,c$.

Most of our original equations needing conversion to resultant form will only have a single 
coefficient $> 1$ and the other two coefficients will = 1.  
So for ease of reference, when only a single coefficient $> 1$, that coefficient will be mapped to 
coefficient $f$ (itself associated with base $c$ to the $z$-th power) and coefficients $d$ and $e$ 
will be = 1. 
Each of $a,b,c$ will always have a coefficient $\geq$ 1, so when a coefficient of an original equation 
is referred to as \emph{the} coefficient, it specifically means 
the 
sole 
coefficient $> 1$.

Regarding the multiplier $N$, consider the case that has only a single coefficient:
$\pm a^x \pm b^y \pm f c^z = 0$
.
Using the symmetrical form of the equation allows us to avoid repeating some arguments and 
also allows us to consider only positive integers throughout this paper.
In order to convert this original equation into its resultant form, the equation must be 
multiplied by a number $N$ such that:
$N$ is an $x$-th power, $N$ is a $y$-th power and $N f$ is a $z$-th power.  One method of 
generating such an $N$ is 
let $N = f^q$ such that $f^q$ is an $x$-th power, $f^q$ is a $y$-th power, and $f^{q+1}$ is a 
$z$-th power.  We see that
$q \equiv 0 \pmod{x};$
$q \equiv 0 \pmod{y};$
$q + 1 \equiv 0 \pmod{z}\enspace.$
The Chinese Remainder Theorem (CRT) guarantees a solution for $q$ with this system of congruences 
when $\gcd(x y,z) = 1$.
As an example, consider the original equation
$\pm a^x \pm b^y \pm f c^z = 0$ with
$x=y=3;$
$z=5.$
Let $N = f^q$.
Solving the congruences
$q \equiv 0 \pmod{x};$
$q \equiv 0 \pmod{y};$
$q \equiv -1 \pmod{z}$
per the CRT gives $q = 9$, and multiplying the equation by $N=f^9$ yields a conversion 
from the original to the resultant form:
$(f^3 a)^3 + (f^3 b)^3 = (f^2 c)^5\enspace.$

However, an $N$ produced in this fashion is not necessarily the smallest number possible that meets 
those conditions. Consider the situation when the coefficient happens to be a perfect square (such as: 
$1369 * 39^5 + 22505^3 = 22586^3$) with $g^2=f$,
$\pm a^3 \pm b^3 \pm g^2 c^5 = 0\enspace.$
Instead of multiplying by the entire coefficient to the ninth power $(g^2)^9$, 
we can instead multiply by $g^3$, giving
$\pm (ga)^3 \pm (gb)^3 \pm (gc)^5 = 0\enspace.$
This multiplier reduction applies not only when a coefficient is a perfect power, but also to any perfect 
power that divides the coefficient.  To produce the smallest $N$ such that $N$ is a perfect 
$x$-th power, $N$ is a perfect $y$-th power, 
and $N f$ is a perfect $z$-th power, we need to consider the highest power of each distinct prime factor of $f$.

Let $p > 1$ be a prime number.  For an integer $n \geq 1$, the \emph{p-adic valuation} $v_{p}$($n$) 
is defined to be the largest integer $r \geq 0$ such that $p^r \mid n$.
Consider an equation with 
$d \geq 1;$
$e \geq 1;$
$f > 1$.  
To convert from its original form to 
the resultant form, for every distinct prime $p \mid f$, the entire equation 
must be multiplied by number $p^q$ such that
$q \equiv $-$v_{p}$($d$) $\pmod{x};$
$q \equiv $-$v_{p}$($e$) $\pmod{y};$
$q \equiv $-$v_{p}$($f$) $\pmod{z}\enspace.$
But $d,e,f$ are pairwise coprime, so any prime dividing $f$ will not divide $d$ or $e$, so the congruence reduces to:
$q \equiv 0 \pmod{x};$
$q \equiv 0 \pmod{y};$
$q \equiv $-$v_{p}$($f$) $\pmod{z}\enspace.$
When $\gcd(x y,z)=1,$ the CRT guarantees a minimal unique solution.  When $\gcd(x y,z) \neq 1,$ 
and $\gcd(d,e,f) = 1,$ the CRT guarantees that no solution exists.

To complete the conversion of the equation from its original form to the resultant form, the 
coefficients $d$ and $e$ will also need to be ``spread'' across the resultant equation in a similar 
manner in order to make $A,B,C$ perfect powers.
\\
To convert an equation from its original form
$d a^x + e b^y = f c^z$
to its resultant form, the entire equation must be multiplied by a number $N$, the smallest of which is:
\\
\[
N = 
\prod_{p \mid d} p^{\smallest_{qd}(p)} 
\prod_{p \mid e} p^{\smallest_{qe}(p)} 
\prod_{p \mid f} p^{\smallest_{qf}(p)} 
\]
with the functions:

$\smallest_{qd}(p) =$ smallest $q$ where $q \equiv 0 \pmod {\lcm(y,z)}, q \equiv $-$v_{p}$($d$) $\pmod{x}$

$\smallest_{qe}(p) =$ smallest $q$ where $q \equiv 0 \pmod {\lcm(x,z)}, q \equiv $-$v_{p}$($e$) $\pmod{y}$

$\smallest_{qf}(p) =$ smallest $q$ where $q \equiv 0 \pmod {\lcm(x,y)}, q \equiv $-$v_{p}$($f$) $\pmod{z}$
\\
\\
If these systems of congruences each have a solution, the CRT can guarantee a minimum product $N$.
Table~\ref{table1} shows the minimum multiplier information obtained by the CRT for all
exponent sets with each exponent $\leq$ 5.

\begin{table}
\caption{
Minimum prime power multiplier for each prime dividing the coefficient.
Given the exponent set \{$x,y,z$\} shown in column~1, for the original equation
$\pm d a^x \pm e b^y \pm f c^z = 0$,
for each $p$ dividing the coefficient shown in column~2, 
this table shows the minimum $q$ for $p^q$ that the equation needs to be multiplied by 
to convert the original equation into a resultant equation.
\label{table1}}
\centering
\begin{tabular}{ c c r r r r }
\hline\noalign{\smallskip}
Exponent & Coeff associated with   & \multicolumn{4}{ c }{$v_{p}$(coefficient)}\\
set      & base to this exponent & 1 & 2 & 3 & 4\\
\noalign{\smallskip}
\hline
\{4,4,3\} & 3 & 8 & 4 & \cellcolor[gray]{0.5} & \cellcolor[gray]{0.5} \\
\hline
\{5,5,3\} & 3 & 5 & 10 & \cellcolor[gray]{0.5} & \cellcolor[gray]{0.5} \\
\hline
\{3,3,4\} & 4 & 3 & 6 & 9 & \cellcolor[gray]{0.5} \\
\hline
\{5,5,4\} & 4 & 15 & 10 & 5 & \cellcolor[gray]{0.5} \\
\hline
\{3,3,5\} & 5 & 9 & 3 & 12 & 6 \\
\hline
\{4,4,5\} & 5 & 4 & 8 & 12 & 16 \\
\hline
\{3,4,5\} & 3 & 20 & 40 & \cellcolor[gray]{0.5} & \cellcolor[gray]{0.5} \\
\hline
\{3,4,5\} & 4 & 15 & 30 & 45 & \cellcolor[gray]{0.5} \\
\hline
\{3,4,5\} & 5 & 24 & 48 & 12 & 36 \\
\hline
\end{tabular}
\end{table}

\subsection{Non-coprime Exponents}
\label{subsection1dot-coprimeexpo}

The CRT guarantees that an $N$ exists to convert a given original equation when the moduli are pairwise
coprime (in our case, when the exponents are pairwise coprime).  But when exponents are not pairwise
coprime, $N$ can only exist when the coefficients of the bases of the non-pairwise coprime exponents
have factors in common.  But being that $\gcd(d,e,f) = 1$, we know that:
\begin{compactitem}
\item{
If an original equation has two or more bases with pairwise non-coprime exponents, the equation can only 
be converted into a resultant form when the coefficients of each of those bases $= 1$.
}
\item{
Correspondingly, the original equation cannot be converted if all three exponents have a 
factor in common (because then all coefficients $= 1$, and we know that at least one of 
$d,e,f$ must be $> 1$).
}
\end{compactitem}

\subsection{$x$-th, $y$-th, $z$-th Power Free Coefficients}
\label{subsection1dot-powerfree}

An integer $n$ is said to be \emph{k-free} ($k \geq 2$) if for every prime $p$ the p-adic valuation 
$v_{p}$($n$) $< k$ (that is, $p^k \nmid n$).
We can specify that original equation coefficients are minimized and their associated bases are 
maximized such that:
$d$ is $x$-th power free
,
$e$ is $y$-th power free
,
$f$ is $z$-th power free.
Having an equation with
$v_{p}$($d$) $\geq x$, $v_{p}$($e$) $\geq y$, or $v_{p}$($f$) $\geq z$ and converting 
it to a resultant equation yields a final equation that is identical to that produced 
by first moving powers from the coefficients to their associated bases such that
$v_{p}$($d$) $< x$, $v_{p}$($e$) $< y$, $v_{p}$($f$) $< z$, then multiplying that 
equation by a number that is a perfect $x$-th power, perfect $y$-th power, and perfect $z$-th power.
Consider this example:
Exponent set \{$x,x,z$\} with
$\gcd(x,z) = 1;$
$d=e=1;$
and coefficient $f$
is the product of non-powered prime: $p_{1}$, and 
contains a $z$-th powered prime: $p_{2}$: $f = p_{1} p_{2}^z$.
As part of the process of converting this original equation into its resultant form, 
the entire equation must be multiplied by $p_{2}^q$ such that $p_{2}^q$ is an $x$-th power 
and $p_{2}^{q+z}$ is a $z$-th power.  So:
$q \equiv  0 \pmod{x};$
$q \equiv -z \pmod{z}\enspace.$
This is equivalent to looking for:
$q \equiv  0 \pmod{x};$
$q \equiv  0 \pmod{z}\enspace.$
With $x$ and $z$ coprime, the smallest $q$ is $x z$. With the allowance that the
entire resultant equation can always be further multiplied by a number, 
there is no reason to allow coefficients to have a $v_{p}$(coefficient) $\geq$ the 
exponent of the associated base.
\\
\\
Putting all this information together, in order for 
an Original Pegg Equation to be converted to a Resultant Pegg Equation with a Pegg Value $> 1$,
it must be of the form:

$d a^x + e b^y = f c^z$
\qquad
\enspace
\thinspace
\thinspace
$x,y,z \geq 3$
\qquad
\qquad
\qquad
\qquad
$a,b,c,d,e,f \geq 1$

$\gcd(d a, e b, f c) = 1$
\qquad
at least one of $d,e,f > 1$
\qquad
$\gcd(x,y,z) = 1$

$d$ is $x$-th power free; 
$e$ is $y$-th power free; 
$f$ is $z$-th power free;
and any two bases that have non coprime exponents must each have an associated coefficient = 1.
\\
Continuing, we can also determine the types of original equations that could be converted to 
resultant equations with a minimum desired Pegg Value.

\subsection{Smallest Possible $\gcd(A,B,C)$}
\label{subsection1dot-smallestgcd}

$\gcd(A,B,C)$ can be no smaller than the resultant coefficient of the base(s) to the 
highest exponent.
Consider an equation with a coefficient $f$ associated with the base to the $z$-th power.
For every prime $p \mid f$, the equation must be multiplied by $p^q$ such that
$p^q$ is an $x$-th power, 
$p^q$ is an $y$-th power, 
and $p^{q+v_{p}(f)}$ is a $z$-th power.
The $p^q$ multiplication will present itself as part
of $D$ (the resultant coefficient as part of the $A$ base) as $p^{q/x}$,
of $E$ (the resultant coefficient as part of the $B$ base) as $p^{q/y}$,
and of $F$ (the resultant coefficient as part of the $C$ base) as $p^{(q+v_{p}(f))/z}$.
Per \refsec\ref{subsection1dot-coprimeexpo}, with $f>1$, we know that $\gcd(z,x) = 1,$ and therefore $z \neq x$.
Consider the two possibilities:
\begin{compactitem}
\item{
With $z < x$:  Because $z < x$, $q/z$ always $> q/x$, so when $q/x$ is an integer, $q/z + v_{p}(f)/z$ must be a greater integer.
}
\item{
With $z > x$: Because $z > x$, $q/z$ always $< q/x$, so when $q/x$ is an integer, $q/z$ must be less than that integer.
The sum $q/z + v_{p}(f)/z$ may $= q/x$, but it cannot $> q/x$, because in order for the sum to equal a higher integer,
$v_{p}(f)/z$ must be $> 1$, and that is impossible (as
$f$ is $z$-th power free, $v_{p}(f) < z$, and $v_{p}(f)/z < 1$).
}
\end{compactitem}
The above logic applies to all combinations of bases and exponent sets and primes dividing each coefficient.
An original equation must be one of the following exponent sets:
\begin{compactitem}
\item{
a single base has the highest exponent, and the base to the highest exponent must have a coefficient (such as \{3,3,4\}).
}
\item{
a single base has the highest exponent, and it may or may not have an associated coefficient (such as \{3,4,5\})
}.
\item{
two bases have the highest exponent, and each of those bases are precluded from having a coefficient (such as \{5,5,3\}).
}
\end{compactitem}
Regardless of the situation, the resultant coefficient of the base(s) associated with the the highest power will
divide the resultant coefficients of the remaining base(s).

\subsection{Re-associating $\min(A,B,C)$}
\label{subsection1dot-move}

If $\min(A,B,C)$ is not associated with a base to the highest exponent, the entire equation can be 
multiplied by a number that will make $\min(A,B,C)$ associated with a base to the highest exponent.  
For example, consider the example Original Pegg Equation

$5 ^ 3 + 427 ^ 3 = 60073 * 6 ^ 4$.
\\
It converts to the smallest resultant equation

$(60073 * 5)^3 + (60073 * 427)^3 = (60073 * 6)^4$
\\
which does not have $\min(A,B,C)$ associated with the base to the highest power.  This equation 
has a Pegg Value $= a = 5$.  Further multiplying the resultant equation by $N = 2^{\lcm(x,y,z)} = 2^{12}$ yields

$(2^4 * 60073 * 5)^3 + (2^4 * 60073 * 427)^3 = (2^3 * 60073 * 6)^4$
\\
which moves $\min(A,B,C)$ to $C$, but the Pegg Value $<c$ as the gcd ``steals'' a factor of 2 from the $c$ base,
resulting in a Pegg Value $=3$.  
So we do want to multiply by a number that will not ``steal'' a factor from the original base.
If the resultant equation was instead multiplied by $N = 5^{12}$, 
the final equation would be

$(5^4 * 60073 * 5)^3 + (5^4 * 60073 * 427)^3 = (5^3 * 60073 * 6)^4$
\\
and the Pegg Value $= c = 6$ (at the cost of a much larger resultant equation).

\subsection{Reduced Pegg Value When Converting from Original to Resultant Form}
\label{subsection1dot-steal}

This ``stealing'' can also result when the equation is multiplied by the $N$ to convert from its 
original form to resultant form.  
It can occur any time the $p^q$ as part of the resultant coefficient
of the base to the highest power is of a smaller power than the $p^q$ as part of the
resultant coefficient of both the other bases.
As an example, 
Table~\ref{table2} shows that this can occur with the \{3,3,5\} exponent set
with $v_{p}(f)$ either 1 or 3.
Consider the original \{3,3,5\} equation:

$(5 * 23)^3 + (2^7)^3 = f (3)^5$, with $f=3 * 7 * 709.$
\\
After converting it to its resultant form:

$(f^3 * 5 * 23)^3 + (f^3 * 2^7)^3 = (f^2 * 3)^5$
\\
the resulting $\min(A,B,C) = C$, but instead of $\gcd(A,B,C)$ being $f^2$, the 
gcd is $f^2 3$, because it ``steals'' a factor of 3 from $c$.  So the Pegg Value 
of the resultant equation is 1, less than the smallest original base to 
the highest exponent.
And there is no away around this ``stealing'' during the original to resultant 
equation conversion, because we must multiply the equation by the primes that 
divide the coefficient.

\begin{table}
\caption{
Power of prime in resultant coefficient after conversion from original equation to resultant equation.
Given the exponent set \{$x,y,z$\} shown in column~1, for the original equation
$\pm d a^x \pm e b^y \pm f c^z = 0$,
for each $p$ dividing the coefficient 
shown in column~2, this table shows how the multiplier $p^q$ (with $q$ determined from
Table~\ref{table1}) is represented in the resultant coefficients
using the format [log$_p(D)$, log$_p(E)$, log$_p(F)$].
\label{table2}
}
\centering
\begin{tabular}{ c c c c c c c }
\hline\noalign{\smallskip}
Exponent & Coeff associated with   & \multicolumn{4}{ c }{$v_{p}$(coefficient)} \\
set      & base to this exponent & 1 & 2 & 3 & 4 \\
\noalign{\smallskip}
\hline
\{4,4,3\} & 3 &   [ 2, 2, 3] &   [ 1, 1, 2] & \cellcolor[gray]{0.5} & \cellcolor[gray]{0.5} \\
\hline
\{5,5,3\} & 3 &   [ 1, 1, 2] &   [ 2, 2, 4] & \cellcolor[gray]{0.5} & \cellcolor[gray]{0.5} \\
\hline
\{3,3,4\} & 4 &   [ 1, 1, 1] &   [ 2, 2, 2] &   [ 3, 3, 3] & \cellcolor[gray]{0.5} \\
\hline
\{5,5,4\} & 4 &   [ 3, 3, 4] &   [ 2, 2, 3] &   [ 1, 1, 2] & \cellcolor[gray]{0.5} \\
\hline
\{3,3,5\} & 5 &   [ 3, 3, 2] &   [ 1, 1, 1] &   [ 4, 4, 3] &   [ 2, 2, 2] \\
\hline
\{4,4,5\} & 5 &   [ 1, 1, 1] &   [ 2, 2, 2] &   [ 3, 3, 3] &   [ 4, 4, 4] \\
\hline
\{3,4,5\} & 3 &   [ 7, 5, 4] &   [14,10, 8] & \cellcolor[gray]{0.5} & \cellcolor[gray]{0.5} \\
\hline
\{3,4,5\} & 4 &   [ 5, 4, 3] &   [10, 8, 6] &   [15,12, 9] & \cellcolor[gray]{0.5} \\
\hline
\{3,4,5\} & 5 &   [ 8, 6, 5] &   [16,12,10] &   [ 4, 3, 3] &   [12, 9, 8] \\
\hline
\end{tabular}
\end{table}

\subsection{Minimum Original Equation Bases for a Given Pegg Value}
\label{subsection1dot-smallest-bases}

Per \refsec\ref{subsection1dot-smallestgcd}, $\gcd(A,B,C)$ can be no smaller than the resultant 
coefficient associated with the base(s) to the highest power.  And per \refsec\ref{subsection1dot-move},
$\min(A,B,C)$ can always be associated with the smallest base to the highest power.  
Therefore the smallest base
criteria for an original equation to convert to a resultant equation of a given Pegg Value $V$ is only that
each of the base(s) to the highest exponent must be $\geq V$
.
Under the condition that the equation not be multiplied by a number to re-associate $\min(A,B,C)$ to 
a different base, $\min(A,B,C) / \gcd(A,B,C)$ must immediately be $\geq V$.
Consider the original \{3,4,5\} equation with $x=3, y=4, z=5$ and a coefficient of 2 associated with the base
to the third power: $\pm 2 a^3 \pm b^4 \pm c^5 = 0.$  It
converts to its resultant form:
$\pm (2^7 a)^3 \pm (2^5 b)^4 \pm (2^4 c)^5 = 0.$
For this equation to be able to be converted to a resultant equation with a Pegg Value $\geq V$,
without further multiplying the equation to re-associate $\min(A,B,C),$
the following conditions must be met:

$D a / F \geq V$ 
\enspace
\enspace
\enspace
\thinspace
$2^7 a / 2^4 \geq V$ 
\qquad 
$a \geq V / 8$

$E b / F \geq V$ 
\qquad 
$2^5 b / 2^4 \geq V$ 
\qquad 
$b \geq V / 2$

$F c / F \geq V$ 
\qquad 
$2^4 c / 2^4 \geq V$ 
\qquad 
$c \geq V$

The Pegg Value may still be less than $V$ due to gcd ``stealing'' when $\gcd(2,c) \neq 1$, but
the above minimums must be met for the equation to have a Pegg Value $\geq V$.

\section{Generating Equations with a Desired Pegg Value}
\label{section4}

Darmon and Granville \cite{bib-darmon-granville} showed that with fixed coefficients and exponents,
there can be, at most, finitely many Original Pegg Equations with unknown integers $a,b,c$.
Oesterlé and Masser's ABC-conjecture 
implies that for fixed $d,e,f$ coefficients, even allowing $x,y,z$ to vary, 
the total number of solutions is limited, thereby implying that the total number of T-Z Conjecture counterexamples 
is finite.

But for our Pegg Value searching, the sizes of the coefficients on an original equation do not have 
an impact on its Pegg Value - only its resultant size.  
As $d,e,f$ are not fixed in Original Pegg Equations, the limitation proven by Darmon and Granville 
does not constrain the number of 
original equations.

With $W = V^{x+2}-1$, consider the identity
$(W^{x+2})^{x} + (W^{x+1})^{x+1} = (W^{x} V)^{x+2}$.
As $W$ and $V$ are necessarily coprime
\footnote{
Note that with $J + K = L$, any two terms being coprime forces all three to be pairwise coprime.
$W$ and $V$ are coprime, because they are both coprime to 1, the third term in the equation $W + 1 = V^{x+2}$
}
, $\gcd(A,B,C) = W^x$, and the Pegg Value = $V$ whenever 
$W^{x} V < W^{x+1}\enspace.$
With $x \geq 3$, this is true whenever $V > 1$.
So with $V \geq 2$, $x \geq 3$, the equation
yields a Resultant Pegg Equation with exponents 
\{$x,x+1,x+2$\} and a Pegg Value of $V$.
Therefore the answer to Pegg's question:
``What is the maximal value of $\min(A,B,C) / \gcd(A,B,C)?$''
is: there is no maximal value; we can construct an equation with any desired Pegg Value.
Using the above method with $x = 3$, we can generate a resultant equation with a Pegg Value
of 60000 and a size $\approx 2^{1270}$.  This is far from proving that 60000 is the highest
Pegg Value for equations $\leq 2^{1270}$.
Indeed, we can use a different identity on the \{3,3,5\} exponent set to generate smaller equations
of a given Pegg Value.
The obvious question arises:  What is the highest 
Pegg Value for equations with $C^z \leq$ a certain size?
\section{The Highest Pegg Value in All Equations $\leq 2^{100}$}
\label{section5}

The strategy is to guess the exponent set that will generate the highest Pegg Value within a given 
equation range.  Then search the range for equations with higher and higher Pegg Values.  Each 
time a new solution is found, it reduces the search space.  For example, after finding an equation 
with a Pegg Value of 1000, we can exclude from the search any original equations that could not possibly 
convert to a resultant equation with a Pegg Value $> 1000$.

After the maximum Pegg Value is found for the given exponent set within the resultant equation size 
range, this Pegg Value can be used to reduce (or eliminate) the search space of other exponent sets.  
A quick guess as to the type of exponent sets of original equations that will yield the highest 
Pegg Value:

\begin{compactitem}
\item{
Will have small exponents.
In addition to the ABC-conjecture's expectation that possible solutions thin out as 
exponents become higher, within a given equation range higher exponents reduce the 
size of their respective base, thereby limiting the equation's Pegg Value.
}

\item{
Will have two identical exponents.
We know that the exponent set \{3,4,5\} yields an infinite number of Pegg Value solutions, 
but any coefficient on this set will necessarily need a large multiplier to make the coefficient 
become a perfect 3rd, 4th, and 5th power.  This large multiplier significantly reduces the 
maximum Pegg Value within a given range.
}

\item{
Will have an exponent set where the multiplier needed to convert between 
the original equation and its resultant equation will be as small as possible.
For example, exponent set \{3,3,4\} requires a square-free coefficient to be cubed to 
convert an original equation to resultant form.  But exponent set \{4,4,3\} requires
a square-free coefficient to be 8th powered to convert an original equation to resultant form.
}
\end{compactitem}

These conditions indicate that exponent set \{3,3,4\} is the best candidate exponent set.  
To gain some confidence in this assessment, we searched the abc@home database of 7.5 million ABC-Hits \cite{bib-abc-hit-database}
looking for hits that could be formed into Original Pegg Equations 
(with each exponent $\leq 5$)
that had a large Pegg Value in relationship to their resultant equation size.

Recall that if non-zero positive integers $A + B = C, A < B < C$, 
and gcd($A,B,C$) $= 1$, then the three are called an \emph{ABC-Triple}.
The \emph{radical} of $N$ is the product of the distinct primes dividing $N$.
i.e., the largest square-free factor of $N$.  
The \emph{ABC-Power} of an ABC-Triple is defined as log(C) / log(rad($A B C$)).
If the ABC-Power of an ABC-Triple is $> 1$, then the triple is called an \emph{ABC-Hit}.  

Original Pegg Equations are not necessarily ABC-Hits. For example, 
$61^3 + 67^3 = 4123 * 2^7$
converts to a resultant equation with a Pegg Value of 2, yet has an ABC-Power of only 0.7602.
But ABC-Hits 
are excellent candidates to have a large 
Pegg Value in relationship to their resultant equation size (because original equations with small radicals 
will tend to be those that have small coefficients).

We define the \emph{Pegg Power} to be 
log(Pegg Value of equation) / log(resultant equation size).  
Similar in concept to the ABC-Power, the Pegg Power provides an easy guide as to the Pegg Value 
``quality'' of the equation.
The highest Pegg Power found was for the
\{3,3,4\} 
Original Pegg Equation
$14 * 111 ^ 4 + 3595 ^ 3 = 3649 ^ 3 $
which converts to a resultant equation
with a Pegg Value of 111 and a Pegg Power = 0.1448.
This resultant equation had the highest Pegg Power of all ABC-Hits in the database 
and is shown in Table~\ref{table-best334} as it has the highest Pegg Value of all \{3,3,4\} equations $< 3558^3$.
Forty-nine of the top fifty Pegg Powers in the database were \{3,3,4\} equations.
The highest non-\{3,3,4\} Pegg Power in the database was thirteenth place 0.10677 for a \{3,3,5\} equation.

This information lends support to the plan 
to find the highest Pegg Value in \{3,3,4\} equations 
$\leq 2^{100}$, and use that result to reduce or eliminate the search space in other exponent sets.

\subsection{Perfect Power Testing}
\label{subsection5dot-powertest}

In order to search for the highest Pegg Value equation $\leq 2^{100}$, we will need a perfect power 
tester for cubes, fourth powers, and fifth powers.
But each power testing algorithm can be distinct.  Having three different perfect power testers reduces 
the runtime inefficiency of having a
general-purpose perfect power tester that checks an input for being a perfect $y$-th power against multiple $y$.

Our perfect power testers work along the lines used by GMP for testing perfect squares \cite{bib-gmp-perf-square}, 
namely, before performing a rigorous but expensive $y$-th power test on an input, first verify 
that the input is compatible with being a $y$-th power modulo small integers.  Most non-powers 
will quickly be eliminated, and for those inputs that pass the residue testing, the algorithm 
will then perform a rigorous perfect power test.

For perfect cube testing, primes 
$p \equiv 1 \pmod{ 3}$
each rule out
$\frac{2}{3} (p-1) / p$
possible inputs.
In addition, powers of an unused prime can assist.  For cubic residues,
the modulus 9 rules out 6/9 of possible inputs.  A modulus of 27 rules out 20/27 of possible inputs,
but being that $\gcd(9,27) \neq 1$, the eliminations overlap,
so it is only a minimal advantage to use 27 as a modulus when also using a modulus of 9.
Using a somewhat arbitrary cut-off point,
checking inputs for being compatible with being a cube modulo 9
and the 34 primes
$p \equiv 1 \pmod{ 3}$
with $p \leq 367$ rules out 99.99999999999999446\% of all inputs.

Notice that 4 is the only exponent $\geq 3$ that needs to be checked that is not a prime.  
Checking for fourth powers against residues for primes
$p \equiv 1 \pmod{ 2} ($and not$ \equiv 1 \pmod{ 4})$
would eliminate
$\frac{1}{2} (p-1)/p$
candidates per prime, but we can do better by only checking against residues for primes
$p \equiv 1 \pmod{ 4}$
 which eliminate
$\frac{3}{4} (p-1)/p$
candidates per prime.  Checking inputs for being compatible with being a fourth power modulo 
9, 16, 49, and the 25 primes
$p \equiv 1 \pmod{ 4} \leq 257$
rules out 99.99999999999999516\% of all inputs.

For fifth powers, checking inputs for being compatible with being a fifth power modulo 25,
and the 23 primes
$p \equiv 1 \pmod{ 5} \leq 521$
rules out 99.99999999999999571\% of all inputs.

Residue testing in such a manner offers another benefit.
Because all our perfect power testing will ultimately be done on inputs of the
form 
$a^x - f c^z$
and
\mbox{$f c^z - a^x$,}
we do not need to compute the actual difference before performing residue
tests.  We can use modular powering, and only compute the true total when we need the rigorous 
perfect power check. For more runtime efficiency, all the modular powering operations can be replaced with 
precalculated table lookups so no actual modular powerings are needed.

Our computer algorithm would instinctively seem to require expensive division instructions to 
perform residue calculations, but divisions needed for the 
modulus operation (\% in C) can be avoided entirely by structuring the code to allow the compiler 
to optimize away actual divisions and instead use multiplicative inverses \cite{bib-amd,bib-coetzee}.  
In GCC, this requires using a hard-coded literal as a modulus (not stored in an array -- even 
a constant array).

The smallest allowable Pegg Equation exponent is 3, so using bases up to $2^{64}-1$ allows us to handle all 
inputs through $(2^{64}-1)^3$ while still using the natively supported
64-bit unsigned integer size of the host computer.  This also contributes to impressively minimal testing time.

As we can perform the residue tests with multiplicative inverses instead of divisions, and can perform memory 
lookups instead of modular powering operations, the time to 
perform a perfect $y$-th power test for non-$y$-th power inputs is effectively independent of 
the size of the inputs $a,x,f,c,z$, and $y$.

\subsection{Algorithm to Find Resultant Pegg Equations with a Minimum Desired Pegg Value}
\label{subsection5dot-Pegg-algo}

This section will present an algorithm that will solve the simple case when two exponents are identical, 
and then conclude with a discussion of how it also works for the more complicated case when all
exponents are different (but still have only a single coefficient).

When two exponents are the same we will map them to the $x$ and $y$ exponents,
allowing us
to follow our convention that 
exponent $z$ is associated with the base with the coefficient.
The equation
$\pm a^x \pm b^y \pm f c^z = 0$ with $x,y,z \geq 3; a,b,c \geq 1$
\\has three permutations when the signs are all positive.  A convenient naming system is to title 
these permutations according to what $b^y$ equals:
\begin{center}
\begin{tabular}{ c c }
\hline\noalign{\smallskip}
Equation Permutation & Permutation Name \\ 
\noalign{\smallskip}
\hline 
$a^x - f c^z = b^y$ & ax\_minus\_cz \\ 
\hline 
$f c^z - a^x = b^y$ & cz\_minus\_ax \\
\hline 
$a^x + f c^z = b^y$ & ax\_plus\_cz \\
\hline 
\end{tabular} 
\end{center}

When two exponents are the same, 
the first and last permutations listed above are identical, and the set of permutations to check 
can be reduced to ax\_minus\_cz and cz\_minus\_ax.

Our algorithmic strategy is to loop through the valid $f$ coefficients, loop through the valid $c$ bases,
then loop through the valid $a$ bases then, depending on the permutation that is being checked, see if
$f c^z - a^x$ or $a^x - f c^z$ yields a perfect $y$-th power.
Consider the exponent set \{3,3,4\} with $x=3$, $y=3$, $z=4$.
There are fewer fourth powers within a given range than there are perfect cubes.  And our perfect 
cube test takes about the same amount of time to perform as our perfect fourth-power test.  So for 
efficiency, our first base loop will range through the fourth powers, and then later perform perfect cube testing.

Our perfect cube test also takes about the same amount of time to perform regardless of the 
size of the input.  So regarding the two bases with identical exponents, it makes sense for 
us to loop through the range of cubes that has a smaller number of cubes than the other range.  

Let $t =$ an offset $>$ 0

Let $r_{min} =$ range minimum $ > 0$ to search for a perfect power

Let $r_{max} =$ range maximum $ > r_{min}$ to search for a perfect power
\\
There are fewer cubes between 
$t + r_{min}$ and $t + r_{max}$
than between
$r_{min}$ and $r_{max}\enspace.$
For efficiency, we establish the convention when $x=y$, the 
$a$ base is larger than the $b$ base.  With permutation type ax\_minus\_cz, we have
$a^x - c^z = b^y$
equivalently
$a^x = c^z + b^y$
.
As applied to our application, there are fewer cubes between
$c^z + b^y_{min}$ and $c^z + b^y_{max}$
than there are between
$b^y_{min}$ and $b^y_{max}\enspace.$
So for the permutation type ax\_minus\_cz, in comparison to looping through the $b$ base and
then checking if the sum $f c^z + b^y$ is a perfect cube,
it is more efficient to loop through the $a$ base, then check if
the difference $a^x - f c^z$ is a perfect cube.
With type cz\_minus\_ax, we have
$f c^z - a^x = b^y$
equivalently
$a^x = f c^z - b^y\enspace.$
Similarly, it will be more efficient to loop through the $a$ base, and check the difference $f c^z - a^x$
and see if the result is a perfect cube.

We can extend the use of perfect power residues even further.  
Consider permutation type cz\_minus\_ax.  
For each $c$ base, instead of looping through the entire range of possible $a$ bases and seeing if 
$f c^z$ - $a^x$ is a perfect $y$-th power, we can
instead loop through only the $a$ bases for a given $c$ base that are compatible with 
$f c^z$ - $a^x$ being congruent to a perfect $y$-th power modulo a few small integers.

Consider again the example exponent set \{3,3,4\}, with 
$x=y=3;$
$z=4.$
and permutation type cz\_minus\_ax with 
$f=8;$
$c=103;$
$a=21.$
We know $f c^z - a^x$ cannot 
be a cube because $8 * 103^4 - 21^3 \pmod{ 7} \equiv 2 \pmod{ 7}$, which is 
incompatible with being a perfect cube (a perfect cube is $\equiv 0,1,$ or $6 \pmod{ 7}$).

Since we know beforehand that with a $f c^z$ residue of $2 \pmod{ 7}$ that an $a$ base of 21 is
incompatible with $f c^z - a^x$ being a perfect cube, there is no sense checking an $a$ base of 21 in the first
place.  We can precalculate a large lookup array that,
for a given $f c^z$,
shows the offsets between $a$ bases that are compatible with the difference being a perfect $y$-th power
(modulo the product of several small integers).  These moduli will be smallest moduli
from the list generated for our perfect power tester.  Increasing the number of small moduli in this
a\_base\_skipahead modulus product increases the size of the table (and table precalculation 
time), but decreases the number of possible $a$ bases that must be checked.  This table gets large
quickly, because it must list the deltas between acceptable $a$ bases (per the a\_base\_skipahead modulus)
for every possible $f c^z$ residue.  Basically, the strategy is to multiply together as many of 
the smallest moduli together to produce a single a\_base\_skipahead modulus product such that the
lookup table will still fit within computer RAM (such as within a typical 4 GiB limit).
Then any $a$ bases that are selected for a given $f c^z$ residue will
automatically be a perfect $y$-th power modulo the a\_base\_skipahead modulus product. 
So the perfect power testing routine will then use remaining moduli that were not used in the
a\_base\_skipahead modulus product.

In addition, sometimes we will know immediately upon checking the current $f c^z$ that there are 
no $a$ bases that could generate a $f c^z - a^x$ congruence compatible with the result being a 
perfect $y$-th power.  For example, if $f c^z \equiv 4 \pmod{ 7}$, there is no $a$ base such that 
cubing it results in $f c^z - a^x$ being congruent to a cube $\pmod{ 7}$.  In this case the entire 
a\_base loop can be avoided.  This lookup table is quite small because it is simply a list of 
$f c^z$ residues that can be skipped.

\begin{algorithm}
\caption{Within the specified equation range, 
find an equation with a minimum desired Pegg Value that is based on an Original Pegg Equation with a single coefficient
\label{algorithm-Pegg}}
\begin{algorithmic}[1]

\FOR{permutation in ax\_minus\_cz \TO cz\_minus\_ax}
  \FOR{$c$ base in all possible inside resultant eq range with this permutation}
    \FOR{$f$ coeff in all possible inside resultant eq range with this permutation and $c$ base}
      \IF {a\_base\_elimination table does not rule out all $a$ bases for this $f$ and $c^z$}
        \FOR{$a$ base in all possible inside resultant eq range with this permutation and $c$ base and $f$ coeff
[use the a\_base\_skipahead lookup table to only check the $a$ bases that are compatible 
with the result being congruent to a $y$-th power modulo a\_base\_skipahead modulus]}
\label{abaseloop}
          \IF { (permutation = ax\_minus\_cz \AND $a^x - f c^z$ = perfect $y$-th power) \OR
                (permutation = cz\_minus\_ax \AND $f c^z - a^x$ = perfect $y$-th power) }
            \IF {$\gcd(a,c) = 1$}
\label{checkgcd}
              \STATE Calculate maximum Pegg Value with this original eq within given resultant eq range
\label{calculatepegg}
              \IF {the maximum Pegg Value $\geq$ minimum desired Pegg Value}
                \RETURN{Pegg Value,Resultant Pegg Equation}
              \ENDIF
            \ENDIF
          \ENDIF
        \ENDFOR
      \ENDIF
    \ENDFOR
  \ENDFOR
\ENDFOR
\RETURN{entire range searched - none exist}

\end{algorithmic}
\end{algorithm}

The specifics of generating the a\_base\_elimination and a\_base\_skipahead tables are 
detailed in Appendix~\ref{appendix-tablemod}.  
The key point is that once the tables are precalculated, they
continually eliminate over 99\% of possible inputs without any further calculations.  
A single increment through
the a\_base\_skipahead table eliminates a large swath of possible inputs.

The specifics of how to calculate the range of 
$c$ bases, 
$f$ coefficients, 
and $a$ bases are
detailed in Appendix~\ref{appendix-coeff}.  Of note is 
a significant reduction in search space for the permutation type cz\_minus\_ax:
$f c^z = a^x + b^y$.
When $x=y$, either $a^x$ or $b^y$
must be $\geq f c^z / 2$, and the other must be $\leq f c^z / 2$.  
This reduces the range to check from the standard range:
$a^x = a_{min}^x .. f c^z - b^y_{min}$  [when $x \neq y$]
to
$a^x = f c^z / 2 .. f c^z - b^y_{min}$  [when $x = y$].
With this reduction, it takes much less time to check the cz\_minus\_ax
permutation than the ax\_minus\_cz permutation.

Notice that the gcd condition on line~\ref{checkgcd} is checked after first verifying that the appropriate 
difference is a perfect power.  If the gcd check was performed before checking for a 
perfect power, it would rule out $\approx 29\%$ of possibilities (approximately $(1 - (6 / \pi^2))$
\cite{bib-hardy}) but even though gcd is an efficient algorithm, we can test for perfect cubes 
even more quickly than calculating gcd when most inputs are not perfect cubes.  And being
that the difference is rarely a perfect cube, it is more
efficient to test the gcd condition after first verifying the perfect $y$-th power condition.

A note on line \ref{calculatepegg}:
``Calculate maximum Pegg Value with this original equation within given resultant equation range''

The first step is to multiply the entire equation by smallest $N$ to produce a valid 
Resultant Pegg Equation.  
At this point, the Pegg Value will be $\geq$ the desired Pegg Value unless either:

\begin{compactitem}
\item{
$\min(A,B,C)$ is not associated with the base to the highest exponent.  In this case, the 
Pegg Value may be able to be increased by multiplying the entire equation by a number that moves 
$\min(A,B,C)$ to a different base.  Of course, the new resulting equation still must fit within 
the desired equation range.  See \refsec\ref{subsection1dot-move}.
}
\item{
$\min(A,B,C)$ is associated with a base to the highest exponent, and the gcd ``steals'' a factor 
from the base associated with the original coefficient.  In this case, nothing can increase 
the Pegg Value.  See \refsec\ref{subsection1dot-steal}.
}
\end{compactitem}
\vspace{4 mm}
This algorithm can also be used for exponent sets where all exponents are different by reordering
the input parameters.  The standard input parameters are \{$x,y,z$\} 
where $x$ and $y$ are exponents associated with bases with no coefficients.

$z = $ exponent associated with the base with coefficient (the algorithm will loop through these $c$ bases)

$x = $ exponent associated with the base with the smaller number of elements (the algorithm will loop through these $a$ bases)

$y = $ exponent associated with the base with the larger number of elements (these are detected by performing perfect $y$-th power tests on $f c^z - a^x$ or $a^x - f c^z$)

When $x=y$, the permutation ax\_plus\_cz was ignored, as it was identical to az\_minus\_cz.  But 
when $x \neq y$, the unprogrammed permutation ax\_plus\_cz can be checked by reordering the inputs
to the algorithm.
Consider the exponent set \{3,4,5\} where the coefficient is
associated with the base to the third power.  
We follow our convention that exponent $z$ is associated with the base with the coefficient, so $z=3$.

The following table shows how we check all three permutation for $z=3$ with our current algorithm 
which handles the two permutation types ax\_minus\_cz and cz\_minus\_ax:

\begin{center}
\begin{tabular}{ c c }
\hline\noalign{\smallskip}
Equation Permutation & Handled with input parameter ordering \{$x,y,z$\} \\ 
\noalign{\smallskip}
\hline 
$a^x - f c^z = b^y$ & ax\_minus\_cz with exponent set \{5,4,3\} \\
\hline
$f c^z - a^x = b^y$ & cz\_minus\_ax with exponent set \{4,5,3\} or \{5,4,3\} \\
\hline
$a^x + f c^z = b^y$ & ax\_minus\_cz with exponent set \{4,5,3\} \\
\hline
\end{tabular} 
\end{center}

Notice that the second possibility can be checked with exponent set \{4,5,3\} or exponent set \{5,4,3\}.
Our ax\_minus\_cz and cz\_minus\_ax search algorithms run more efficiently with $x \geq y$, so it 
will be most efficient to check exponent set \{5,4,3\}.

For the third possibility, it is less efficient at runtime to check ax\_minus\_cz with exponent set \{4,5,3\}
than using a new ax\_plus\_cz routine and checking exponent set \{5,4,3\}, but the inefficiency is not
great enough to overcome the overhead of writing a new routine.  This searching does not take very long
because there are far fewer combinations to check when all three exponents are different.

\begin{figure}[H]
  \centering\includegraphics[scale=0.95]{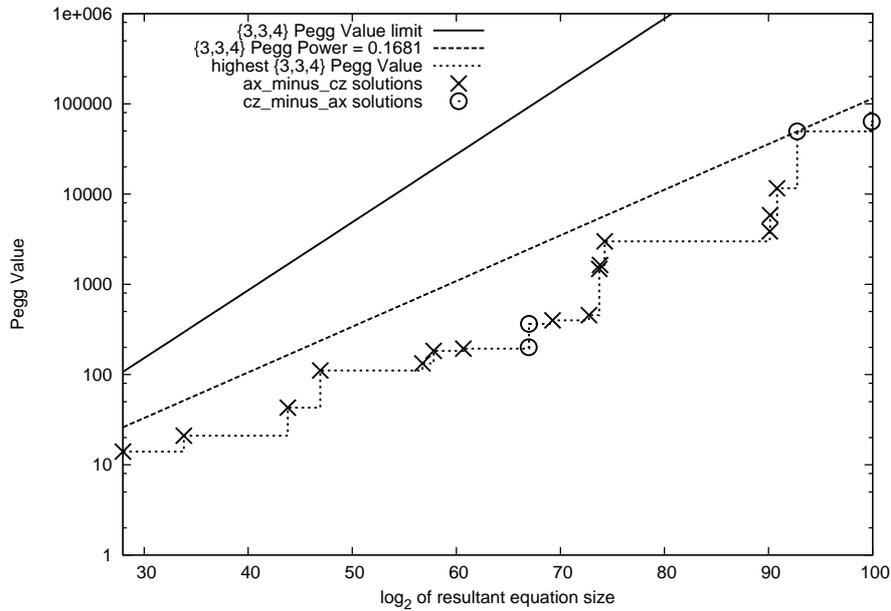}
  \caption{Highest \{3,3,4\} Pegg Values for a given resultant equation size (information from 
Table~\ref{table-best334})
  \label{figure1}}
\end{figure}

Using the above algorithm, it is easy to determine that the smallest \{3,3,4\} 
equation with a Pegg Value $> 1$ is $207 ^ 3 + 126 ^ 4 = 639 ^ 3$
which has Pegg Value of 14.  Searching  \{3,3,4\} equations up to $2^{100}$ for higher and 
higher Pegg Values yields the information in Table~\ref{table-best334}
and Figure~\ref{figure1}.
The last row in Table~\ref{table-best334} shows a resultant equation $\leq 2^{100}$ with a Pegg Value of $63742$.  
The next step is to prove the hunch that this is the highest Pegg Value in all equations 
$\leq 2^{100}$.

\begin{lemma}
There are no T-Z Conjecture counterexamples $\leq 2^{100}$ with a Pegg Value $> 63742$
\label{lemma1}
\end{lemma}

Any equation with a higher Pegg Value must have all bases $> 63742$.  This precludes any base 
having an exponent of 7 or higher, as $63743^7 > 2^{100}$.

A good summary of exponent sets for which it has been proven that no T-Z Conjecture 
counterexamples exist is provided in both \cite{bib-beukers-fermatlectures} and \cite{bib-poonen-twists}.
Of all exponent sets with each exponent $\leq 5$, all but one have been ruled 
out as having any T-Z Conjecture counterexamples.  This information is 
summarized in Table~\ref{table-ruleout}.  Note that exponent 6 is covered under exponent 3, as $a^6 = (a^2)^3$, 
so for T-Z Conjecture counterexamples we only need to consider exponent sets with each 
exponent $\leq 5$.

\begin{table}
\caption{
All combinations of T-Z Conjecture counterexample exponent sets with each 
exponent $\leq 5$, along with a possible reference to a proof showing that none exist
\label{table-ruleout}}
\centering
\begin{tabular}{ c l }
\hline\noalign{\smallskip}
exponent & \multicolumn{1}{c}{reference to proof that exponent set has} \\
set      & \multicolumn{1}{c}{no T-Z Conjecture counterexamples} \\
\noalign{\smallskip}
\hline
\{n,n,n\} & Wiles and Taylor. Fermat's last Theorem. \cite{bib-wiles,bib-taylor} \\
          & although the specific cases with $3 \leq n \leq 5$ were all solved earlier \\
\hline
\{3,3,3\} & Euler (18th century) \\
\hline
\{4,4,4\} & Fermat (17th century) \\
\hline
\{5,5,5\} & Dirichlet and Legendre (19th century) \\
\hline
\hline
\{3,3,4\} & Bruin \cite{bib-bruin-sum-two-cubes} \\
\hline
\{3,3,5\} & Bruin \cite{bib-bruin-sum-two-cubes} \\
\hline
\{4,4,3\} & Lucas (19th century) \\
\hline
\{4,4,5\} & Bruin \cite{bib-bruin-chabauty-methods} showed complete list of \{2,4,5\} \\
\hline
\{5,5,3\} & Poonen \cite{bib-poonen-some-diophantine} \\
\hline
\{5,5,4\} & Poonen \cite{bib-poonen-some-diophantine} ruled out \{5,5,2\} \\
\hline
\hline
\{3,4,5\} & \emph{has not been ruled out} \\
\hline
\end{tabular} 
\end{table}

The only T-Z Conjecture counterexample exponent set listed above that has not been 
ruled out is \{3,4,5\}.  Notice that \{3,4,5\} can be rewritten as a special form of 
\{2,3,5\} for which the complete parameterization was shown by
Johnny Edwards \cite{bib-edwards} and 
conveniently available on Dario Alpern's web-site \cite{bib-alpern}. 
Using these parameterizations, checking all solutions for 

$\pm a^2 \pm b^3 \pm c^5 = 0$ with $a,b,c \geq 1, \gcd(a,b,c) = 1$
\\
shows that for equations $\leq 2^{800}$, in no case is $a$, the base to be squared, itself 
a perfect square.  So there are no \{3,4,5\} T-Z Conjecture counterexamples with an equation 
size $\leq 2^{800}$, thereby proving that there are no T-Z Conjecture counterexamples of any 
exponent set under an equation size $\leq 2^{100}$ that have a Pegg Value $> 63742$.

\begin{table}
\caption{
Smallest \{3,3,4\} resultant equation with a Pegg Value higher than that shown in the preceding row (Row 1 shows the smallest \{3,3,4\} equation with a Pegg Value $> 1$).
\label{table-best334}}
\centering
\begin{tabular}{rrrr}
\hline\noalign{\smallskip}
log$_2$    &       &       &                   \\
(resultant & Pegg  & Pegg  &                   \\
equation)  & Value & Power & Original Equation \\
\noalign{\smallskip}
\hline
27.96 & 14 & 0.1362 &   $23 ^ 3 + 9 * 14 ^ 4 = 71 ^ 3$ \\
\hline
33.81 & 21 & 0.1299 &   $13 * 21 ^ 4 + 163 ^ 3 = 190 ^ 3$ \\
\hline
43.80 & 43 & 0.1239 &   $23 * 43 ^ 4 + 1056 ^ 3 = 1079 ^ 3 $ \\
\hline
46.92 & 111 & 0.1448 &   $14 * 111 ^ 4 + 3595 ^ 3 = 3649 ^ 3 $ \\
\hline
56.75 & 133 & 0.1243 &   $1157 ^ 3 + 139 * 133 ^ 4 = 3558 ^ 3 $ \\
\hline
57.82 & 183 & 0.1300 &   $1966 ^ 3 + 121 * 183 ^ 4 = 5233 ^ 3 $ \\
\hline
60.68 & 194 & 0.1252 &   $126 * 194 ^ 4 + 9071 ^ 3 = 9743 ^ 3 $ \\
\hline
66.96 & 201 & 0.1143 &   $5906 ^ 3 + 8809 ^ 3 = 545 * 201 ^ 4 $ \\
\hline
66.98 & 365 & 0.1271 &   $10973 ^ 3 + 15902 ^ 3 = 301 * 365 ^ 4 $ \\
\hline
69.24 & 399 & 0.1248 &   $12146 ^ 3 + 391 * 399 ^ 4 = 22703 ^ 3 $ \\
\hline
72.75 & 455 & 0.1214 &   $513 * 455 ^ 4 + 33247 ^ 3 = 38872 ^ 3 $ \\
\hline
73.74 & 1482 & 0.1429 &   $1609 ^ 3 + 239 * 1482 ^ 4 = 104857 ^ 3 $ \\
\hline
73.81 & 1638 & 0.1447 &   $97103 ^ 3 + 193 * 1638 ^ 4 = 132095 ^ 3 $ \\
\hline
74.25 & 2994 & 0.1555 &   $104 * 2994 ^ 4 + 226199 ^ 3 = 271127 ^ 3 $ \\
\hline
90.12 & 3858 & 0.1322 &   $25031 ^ 3 + 1570 * 3858 ^ 4 = 703271 ^ 3 $ \\
\hline
90.16 & 5838 & 0.1388 &   $729217 ^ 3 + 971 * 5838 ^ 4 = 1148689 ^ 3 $ \\
\hline
90.82 & 11598 & 0.1487 &   $341 * 11598 ^ 4 + 3662591 ^ 3 = 3809903 ^ 3 $ \\
\hline
92.75 & 49476 & 0.1681 &   $7771657 ^ 3 + 8824055 ^ 3 = 193 * 49476 ^ 4 $ \\
\hline
99.91 & 63742 & 0.1597 &   $2192137 ^ 3 + 20440855 ^ 3 = 518 * 63742 ^ 4 $ \\
\hline
\end{tabular}
\end{table}

\begin{lemma}
There are no Resultant Pegg equations $\leq 2^{100}$ with a Pegg Value $> 63742$ with an exponent of 6
\label{lemma2}
\end{lemma}

In an equation $\leq 2^{100}$, the maximum possible base to a sixth power is $\lfloor 2^{100/6} \rfloor = 104031$.
The only reason to use exponent 6 instead of 3 on a particular base would be if one wanted to force
that particular base to be smaller so that it would be selected as $\min(A,B,C).$
When $\min(A,B,C)$ is associated with an exponent of 6,
$\min(A,B,C) / \gcd(A,B,C)$ can be no higher than $\lfloor 104031 / 2 \rfloor = 52015$,
because $\gcd(A,B,C)$ must be $\geq 2$.
As $52015 < 63743$, there is no reason to consider equations with 
an exponent of 6.

\begin{lemma}
There are no Resultant Pegg equations $\leq 2^{100}$ with a Pegg Value $> 63742$ with two exponents the same
\label{lemma3}
\end{lemma}

The Original Pegg Equation condition that $\gcd(x,y,z) = 1$ precludes $z = x$, and as
all exponents to consider are $< 6$, the following sets must be checked:
\{4,4,3\},
\{5,5,3\},
\{3,3,4\},
\{5,5,4\}.
\{3,3,5\}.
\{4,4,5\}.
Regarding the search for a \{3,3,4\} equation with a coefficient of 2,
Henri Cohen \cite{bib-cohen} showed that the equation:

$\pm a^3 \pm b^3 \pm 2 c^2 = 0$, with $a,b,c \geq 1,\gcd(a,b,c) = 1$
\\
can be parameterized by one of four parameterizations.
Using these parameterizations, checking all solutions with $2 c^2 \leq 2^{240}$ shows that in 
no case is $c$, the base to be squared, itself a perfect square.  So there are no \{3,3,4\} 
original equations $\leq 2^{240}$ with a coefficient of 2, thereby eliminating any need for
our algorithm to consider a coefficient of 2 for exponent \{3,3,4\} in a search for any 
resultant equations $\leq 2^{100}$.

The remaining coefficients for \{3,3,4\} and the other listed exponent sets were checked with 
Algorithm~\ref{algorithm-Pegg} and in no case was an Original Pegg Equation found that could be converted to a 
Resultant Pegg Equation $\leq 2^{100}$ with a Pegg Value $> 63742$.

\begin{lemma}
There are no Resultant Pegg Equations $\leq 2^{100}$ with a Pegg Value $> 63742$ with no two exponents the same
\label{lemma-345}
\end{lemma}

As all exponents to consider are $< 6$, the only the set to check is
\{3,4,5\}.  All the other exponent sets to check had two identical exponents, so we mapped the non-identical 
exponent to $z$ according to our convention that the sole base that had a coefficient $> 1$ 
would be mapped to the $c$ base.  For exponent set \{3,4,5\} all bases could have a coefficient $> 1$, 
so it is not obvious how to map the $x,y,z$ exponents.  To ease confusion, for now we will simply map:
$x = 3;$
$y = 4;$
$z = 5.$
\\
First to note is that there is not room to multiply the entire equation by a number 
in order to increase its Pegg Value by moving $\min(A,B,C)$ from being associated with 
a smaller exponent to a larger exponent:
Per \refsec\ref{subsection1dot-smallestgcd}, the smallest possible $\gcd(A,B,C)=F$, 
so the largest Pegg Value for a \{3,4,5\} equation is $\min(A,B,C) / F$.
The highest value for $C$ is $\lfloor 2^{100/5} \rfloor = 1048576$.  So $\min(A,B,C)$ must be $\leq$ 1048576.
As the Pegg Value is limited to $\min(A,B,C) / F$, the Pegg Value $\leq \lfloor 1048576 / F \rfloor$.
For the Pegg Value to be $\geq$ 63743, $F$ must be $\leq \lfloor 1048576 / 63743 \rfloor \leq 16.$
The smallest number by which the entire equation could be multiplied is $2^{\lcm(x,y,z)} = 2^{60}$, 
which would contribute to a resultant coefficient $F$ of $2^{60/5} = 4096$, 
which is greater than the maximum of 16.

Second to note is that \{3,4,5\} can have only single coefficient.
There must be at least one original coefficient $> 1$, therefore the resulting 
$D,E,F$ will all be $> 1$, and we
also know that $F \leq 16$.  We must answer the question, ``Under what conditions can 
an original \{3,4,5\} equation result in
$1 < F \leq 16$?''

Our search is limited to $F \leq 16 \leq 2^4$, so $v_{p}$($F$) $\leq 4$, for every $p \mid F$.
Looking at Table~\ref{table2} for options where a \{3,4,5\} equation has $F \leq$ 4th power, we have:

\{3,4,5\}, with $v_{p}$($d$) $= 1, F = p^4$
\qquad
\{3,4,5\}, with $v_{p}$($e$) $= 1, F = p^3$

\{3,4,5\}, with $v_{p}$($f$) $= 3, F = p^3$
\\
There cannot be coefficients on more than one base (because there is only room for a 
maximum resultant coefficient of $2^4$ for $F$ and the smallest combination of two 
coprime coefficients would yield a resultant coefficient of $F = p_1^3 p_2^3$ which would
be $> 2^4$, the maximum value of F).  And $p$ can only $=2$, as the smallest possible $F$ is
$p^3$ and cubing the next higher base yields $3^3$ which is greater than the maximum $F$ of 16.
The search for \{3,4,5\} original equations that convert 
to resultant equations $\leq 2^{100}$ with a Pegg Value $> 63742$ is limited to:

\{3,4,5\} with a coefficient of 2 associated with the base to the 3rd power.

\{3,4,5\} with a coefficient of 2 associated with the base to the 4th power.

\{3,4,5\} with a coefficient of 8 associated with the base to the 5th power.
\\
As the equation cannot be multiplied by a number to re-associate $\min(A,B,C),$
\refsec\ref{subsection1dot-smallest-bases} allows us to calculate the original base minimums for
each base (and shows the first \{3,4,5\} case worked as an example).
The three coefficient scenarios were all checked 
with Algorithm~\ref{algorithm-Pegg} and
in no case was an 
Original Pegg Equation found that could be converted to a Resultant Pegg Equation 
$\leq 2^{100}$ with a Pegg Value $> 63742$.

\section*{Conclusion}

The last row of Table~\ref{table-best334} shows a \{3,3,4\} Original Pegg Equation that can be converted
to the equation 
$1135526966 ^ 3 + 10588362890 ^ 3 = 33018356 ^ 4$
that has a Pegg Value of 63742
and a 
$C^z \leq 2^{100}.$
Lemmas \ref{lemma1} through \ref{lemma-345} show that there are no equations 
$\leq 2^{100}$ of any exponent set that have a higher Pegg Value.

Searching and proof verifying took place over several weeks on three AMD 64-bit computers;
from two to four cores - 2.5 Ghz to 3.0 Ghz;
4 GiB to 8 GiB of memory;
running NetBSD 4.0, 4.0.1 or 5.0 \cite{bib-netbsd} with the operating system included GCC 4.1.2 or 4.1.3 \cite{bib-fsf};
and using GMP 4.2.3 or 4.3.1 \cite{bib-gmp}.

The highest Pegg Power of the equations in Table~\ref{table-best334} is 0.1681
.
The Pegg Power limit for \{3,3,4\} equations and \{4,4,3\} equations is 
$\lim_{V \to \infty} \log(V) / \log((2 V)^4) = 0.25$.
Any T-Z Conjecture counterexample or other Pegg Equation exponent set 
will have at least one exponent $\geq 5$, which limits the Pegg Power to 0.2.  So 
the highest Pegg Power for all equations lies somewhere within the range 0.1681 and 0.25.  
Rather than simply looking to increase the lower limit by further searching, a better technique 
to reduce this range could help shed light on reducing or eliminating possible 
T-Z Conjecture counterexamples.  What is the highest Pegg Power?

\section{Acknowledgements}

Without the availability of GMP, the GNU Multiple Precision Arithmetic Library, the work that resulted
in this paper would not have been started.  
I thank the authors and maintainers for such a great product.

\section{Electronic Availability}

The Pegg Value searching code, written in C with calls to GMP, is released under the GNU 
General Public License (GPL) and is available at 

http://www.arthurtownsend.com/peggvalue.htm.

\appendix
\section{Generating the a\_base\_skipahead and a\_base\_elimination Lookup Tables}
\label{appendix-tablemod}

This appendix details the suggested 
a\_base\_elimination moduli and a\_base\_skipahead moduli for all exponent sets with each exponent $\leq 5$.
The use of these tables was outlined in \refsec\ref{subsection5dot-Pegg-algo}:

The a\_base\_elimination table is a list of $f c^z$ residues for which there are no $a$ bases that can
generate a $f c^z$ and $a$ combination that is a $y$-th power.  For these $f c^z$, the entire 
a\_base loop can be avoided.

The a\_base\_skipahead table is a large lookup array showing the offsets between possible $a$ bases for a given
$f c^z$ (offsets between $a$ bases that are not incompatible with the $f c^z$ and $a^x$ combination being a perfect $y$-th power
modulo the a\_base\_skipahead modulus).
Any $a$ bases that are selected will
automatically generate an $f c^z$ and $a^x$ combination that is congruent to a perfect $y$-th power (modulo a\_base\_skipahead modulus),
so the perfect power testing routine will then use $y$-th power testing moduli that were not used in the
a\_base\_skipahead modulus product.

The tables based on these moduli were 
used to generate the information in Table~\ref{table-best334} and also used for verifying 
Lemmas~\ref{lemma3} 
and 
\ref{lemma-345}).

\begin{sidewaystable}
\caption{
Suggested a\_base\_elimination moduli and standard a\_base\_skipahead moduli.
\\
\\
For each exponent set \{$x,y,z$\} shown in column~1 and permutation type shown in column~2,
column~3 shows the a\_base\_elimination modulus to generate the list of $f c^z$ residues that 
result in no possible $a^x$ being able to produce a $f c^z$ and $a^x$ combination that is 
perfect $y$-th power.
\\
\\
Column~4 shows the standard a\_base\_skipahead modulus.  This modulus results in an 
a\_base\_skipahead table less than 4 GiB that can handle all coefficients.
\\
\\
For all possible $c$ bases and for all $f$ coefficients $\leq$ 100000, 
columns~5,~6, and ~7 show how many $a$ bases
are eliminated by the listed lookup table as being able to produce an $f c^z$ and $a^x$ 
combination that could be perfect $y$-th power.
Column~5 shows what percentage would be eliminated solely based on the a\_base\_elimination modulus.
Column~6 shows what percentage would be eliminated solely based on the a\_base\_skipahead modulus.
Column~7 shows what percentage were eliminated based on both tables.
\label{table-moduli}}
\centering
\begin{tabular}{ccrrrrr}
\hline\noalign{\smallskip}
 & & & & \multicolumn{1}{c}{\% of inputs}  & \multicolumn{1}{c}{\% of inputs}  & \multicolumn{1}{c}{\% of inputs}  \\
         &             &                                          & \multicolumn{1}{c}{standard}
& \multicolumn{1}{c}{eliminated by} & \multicolumn{1}{c}{eliminated by} & \multicolumn{1}{c}{eliminated by} \\
exponent & permutation & \multicolumn{1}{c}{a\_base\_elimination} & \multicolumn{1}{c}{a\_base\_skipahead} 
& \multicolumn{1}{c}{elimination}   & \multicolumn{1}{c}{skipahead}     & \multicolumn{1}{c}{combination}   \\
set      & type        & \multicolumn{1}{c}{modulus}              & \multicolumn{1}{c}{modulus}
& \multicolumn{1}{c}{modulus}       & \multicolumn{1}{c}{modulus}       & \multicolumn{1}{c}{of moduli}     \\
\noalign{\smallskip}
\hline
\{3,3,4\} & both          & $7 * 3^2$                                 & 13 * 19 * 31 * 37             &  47.149 &  97.596 &  98.729 \\
\hline
\{3,3,5\} & both          & $7 * 3^2$                                 & 13 * 19 * 31 * 37             &  46.956 &  97.596 &  98.725 \\
\hline
\{4,4,3\} & ax\_minus\_cz & $5 * 2^4 * 17$                            & $3^2 * 13 * 29 * 37$          &  77.398 &  98.063 &  99.562 \\
\hline
\{4,4,3\} & cz\_minus\_ax & $5 * 13 * 2^4 * 17 * 3^3 * 29 * 7^2 * 11^2 * 43^2$ & $3^2 * 13 * 29 * 37$ &  91.228 &  99.160 &  99.853 \\
\hline
\{4,4,5\} & ax\_minus\_cz & $5 * 2^4 * 17$                            & $3^2 * 13 * 29 * 37$          &  69.265 &  98.022 &  99.392 \\
\hline
\{4,4,5\} & cz\_minus\_ax & $5 * 13 * 2^4 * 17 * 3^3 * 29 * 7^2 * 11^2 * 43^2$ & 37 * 41 * 53         &  87.963 &  98.042 &  99.764 \\
\hline
\{5,5,3\} & both          & $11 * 5^2 * 31 * 41 * 61$                 & 61 * 71 * 101                 &  87.392 &  98.915 &  99.830 \\
\hline
\{5,5,4\} & both          & $11 * 5^2 * 31 * 41 * 61$                 & 41 * 61 * 71                  &  87.332 &  98.801 &  99.767 \\
\hline
\{4,3,5\} & both          & 13                                        & $7 * 3^2 * 19 * 31$           &   7.101 &  99.475 &  99.512 \\
\hline
\{5,3,4\} & both          & 31                                        & $7 * 3^2 * 13 * 19$           &  49.951 &  99.773 &  99.886 \\
\hline
\{3,4,5\} & both          & 13                                        & $5 * 3^2 * 2^4 * 17 * 29$     &   7.101 &  99.740 &  99.758 \\
\hline
\{5,4,3\} & both          & 11 * 41                                   & $5 * 3^2 * 13 * 2^4$          &  17.002 &  99.638 &  99.699 \\
\hline
\{3,5,4\} & both          & 31                                        & $11 * 5^2 * 41 * 61$          &  49.951 &  99.993 &  99.997 \\
\hline
\{4,5,3\} & both          & 11 * 41                                   & $5^2 * 31 * 61$               &  17.002 &  99.865 &  99.888 \\
\hline
\end{tabular}
\end{sidewaystable}

\begin{table}
\caption{
Suggested single-coefficient a\_base\_skipahead moduli.
For each exponent set \{$x,y,z$\} shown in column~1 and permutation type shown in column~2,
column~3 shows the single-coefficient a\_base\_skipahead modulus.  This modulus results in 
an a\_base\_skipahead table less than 4 GiB that can handle just a single coefficient.  For 
the two exponent sets that take the longest to search, it is 
worthwhile to check the coefficients that result in the most combinations with their 
own individual a\_base\_skipahead table.
Column~4 lists the coefficients that were checked individually 
to generate the information in Table~\ref{table-best334} and for verifying Lemmas~\ref{lemma3} 
and
\ref{lemma-345}.  The list of coefficients in column~4 is sorted by decreasing number of
combinations the coefficient produces.
\label{table-moduli2}}
\centering
\begin{tabular}{ccrrrr}
\hline\noalign{\smallskip}
         &             & \multicolumn{1}{c}{single-coefficient} & \multicolumn{1}{c}{coefficients} \\
exponent & permutation & \multicolumn{1}{c}{a\_base\_skipahead} & \multicolumn{1}{c}{checked}      \\
set      & type        & \multicolumn{1}{c}{modulus}            & \multicolumn{1}{c}{individually} \\
\noalign{\smallskip}
\hline
\{3,3,4\} & both & 7 * 13 * 19 * 31 * 37    & 
3,4,5,6,7,8,9,10,11,12,13 \\
\hline
\{4,4,3\} & both & $5 * 3^2 * 13 * 29 * 37$ & 4,2,9,25,3,36,49,18,100 \\
\hline
\end{tabular}
\end{table}

The bulk of the time 
of the algorithm will be spent on the coefficients that produce the largest number of
combinations of $c$ bases and $a$ bases within a given resultant equation range.  For \{3,3,4\}, these will
be the smallest coefficients.  The \{3,3,4\} a\_base\_skipahead table for all possible 
coefficients is 2.6 GiB for the modulus 13 * 19 * 31 * 37.  The single coefficient ``3'' provides so many 
$c$ base and $a$ base combinations within a given range that it is worthwhile to generate a \{3,3,4\} 
single-coefficient a\_base\_skipahead table just for the coefficient 3.  For a given a\_base\_skipahead 
modulus, having an a\_base\_skipahead table for only a single
coefficient produces a much smaller table than a table that can handle all coefficients.  And 
because the resulting table is much smaller, we can use a larger modulus 7 * 13 * 19 * 31 * 37 and keep
the table size less than our 4 GiB limit.  Using the additional 7 as part of the modulus
rules out $2 / 3 * (7-1) / 7 \approx 57\%$ more $a$ bases than the original a\_base\_skipahead table.  
Note that these savings do not result in 57\% fewer $a$ bases being checked, because some of 
those now rejected would already have been rejected by the a\_base\_elimination modulus
which also includes the modulus of 7.  But the benefits of a custom individual coefficient
a\_base\_skipahead modulus are sufficiently great that they can be used for a few
coefficients on the \{3,3,4\} exponent set (as shown in Table~\ref{table-moduli2}, column~4). It is not 
worthwhile using a custom a\_base\_skipahead table for every coefficient on the
\{3,3,4\} exponent set, because the higher the coefficient the fewer combinations it 
produces, and the benefit of the 57\% $a$ base reduction must be balanced against the
single-coefficient table precalculation time. 
\\
\\
Some notes on the tables:

Regarding the coefficients tested individually (listed in Table~\ref{table-moduli2}, Column 4):
\begin{compactitem}
\item{
Table~\ref{table1} shows 
that regardless of the powers of the prime factors that divide $f$, for exponent 
set \{3,3,4\} the smallest $N$ possible is always $N=f^3$.
So for \{3,3,4\}, 
the smaller the coefficient, the more combinations
there are to be checked. 
}
\item{
Table~\ref{table1} shows that the multiplier to convert a \{4,4,3\} Original Pegg Equation
to its resultant equation is smallest when $v_p(f)=2$, so for \{4,4,3\}, 
the smaller and the ``more square'' the coefficient, the more combinations
there are to be checked. 
}
\end{compactitem}

Exponent set \{3,3,4\} shows the coefficient 7 being included in the list of coefficients 
checked individually.  For coefficient 7, the listed individual a\_base\_skipahead modulus 
7 * 13 * 19 * 31 * 37 offers no reduction in the $a$ base search space over the smaller 13 * 19 * 31 * 37.  
So for coefficient 7, either the standard a\_base\_skipahead table or single-coefficient 
a\_base\_skipahead table could be used with identical effect, or searching efficiency could 
be improved by using a modulus of 9 * 13 * 19 * 31 * 37.

Exponent set \{3,3,5\} does not have coefficients checked individually because after 
increasing the standard a\_base\_skipahead modulus, even for use with just a single 
coefficient, the resulting \{3,3,5\} single-coefficient a\_base\_skipahead table is 
16 GiB (which exceeds the 4 GiB limit).  Other exponent sets do not take long to check, 
so it is not worth the overhead to precalculate single-coefficient a\_base\_skipahead 
tables for them.

For a given exponent set, the ax\_minus\_cz and cz\_minus\_ax a\_base\_elimination moduli are different only 
when two exponents are even, which when each exponent $\leq 5$ means when two exponents $= 4$, which means
\{4,4,3\} and \{4,4,5\}.

When based upon identical a\_base\_skipahead moduli, the generated a\_base\_skipahead tables for 
ax\_minus\_cz and cz\_minus\_ax permutation types for each exponent set are identical except when $y$ is even, 
which when each exponent $\leq 5$ means when $y = 4$, which means \{4,4,3\}, \{4,4,5\}, \{3,4,5\}, \{5,4,3\}.

The a\_base\_skipahead 
table is overwhelmingly size-sensitive to the original moduli size, so for non-prime moduli typically 
the smallest helpful prime power is used (such as 9 for \{3,3,4\} and \{4,4,3\}), even though 
a higher power (such as 27) eliminates slightly more candidates.  But the a\_base\_elimination 
table is not under such strict size restrictions, so it can easily use the higher power if the higher power 
assists. For \{3,3,4\} the a\_base\_elimination modulus 27 offers no advantage over 9, but for 
\{4,4,3\} cz\_minus\_ax, the a\_base\_elimination modulus 27 does eliminate a few more candidates.

The general plan is to use whatever small moduli as part of the a\_base\_elimination modulus that help
eliminate $a$ bases, then use small moduli coprime to those as part of the a\_base\_skipahead modulus.
This combination of coprime moduli typically eliminates the most candidates.  But some exceptions are noted:
Exponent set \{4,4,3\} permutation type cz\_minus\_ax uses the modulus 29 as part of both the 
a\_base\_elimination modulus and the a\_base\_skipahead modulus. Instead of using 29, if the 
a\_base\_skipahead modulus was to use the next smallest unused perfect $y$-th power modulus (41) 
the resulting a\_base\_skipahead table would be $>$ 4 GiB.
Keeping the modulus of 29 in the a\_base\_elimination modulus provides for quicker elimination 
for those $f c^z$ residues that do not have any valid $a$ base combinations because the 
a\_base\_elimination table is more likely to be stored in cache.  And keeping the modulus of 
29 in the a\_base\_skipahead table allows the
a\_base\_skipahead modulus to reduce the number  of $a$ bases for those $f c^z$ residues that
do have valid $a$ bases.  
Similarly, \{4,4,3\} permutation type cz\_minus\_ax uses $3^2$ as part of the a\_base\_skipahead modulus
and, as mentioned in the previous item, uses $3^3$ as part of the a\_base\_elimination modulus.

The a\_base\_elimination modulus may not exclude as many possibilities as possible.  For example, for 
\{4,4,3\} permutation type cz\_minus\_ax, the a\_base\_elimination 
modulus contains a modulus of 16, even though 1024 would rule out a few more candidates.
The small number of candidates that would be excluded by a modulus of 1024 that are not excluded by 
the modulus of 16 are not worth the overhead of the increased a\_base\_elimination table size.  
Instead, these few candidates that make it past the a\_base\_elimination
check will be immediately rejected by the a\_base\_skipahead table.

Exponent set \{4,4,5\} permutation type cz\_minus\_ax could use the same a\_base\_skipahead modulus 
as \{4,4,5\} permutation type ax\_minus\_cz, but the \{4,4,5\} permutation type cz\_minus\_ax 
a\_base\_skipahead modulus is modified to be the product of moduli that are coprime to the 
exponent set's a\_base\_elimination modulus.

Exponent set \{5,5,4\} could use the same 61 * 71 * 101 a\_base\_skipahead modulus as exponent 
set \{5,5,3\}, but set \{5,5,4\} takes such a short time to check all possible equations that it 
is not worth the extra table precalculation time.

All the a\_base\_skipahead table offsets fit within a 16-bit unsigned table entry (all 
offsets $\leq 65535$), other than some of the \{4,4,3\} single-coefficient 
a\_base\_skipahead tables.  For some \{4,4,3\} single-coefficient tables
(such as for coefficient 3, permutation cz\_minus\_ax), the a\_base\_skipahead 
modulus eliminates so many $a$ base candidates that sometimes the delta between valid $a$ bases 
for a given $f c^z$ residue is $> 65535$, requiring the entire table to use 4-byte fields instead 
of the 2-byte fields required for all the other tables.

As mentioned in \refsec\ref{subsection5dot-Pegg-algo}, the a\_base\_elimination lookup table is quite 
small because it is simply a list of $f c^z$ residues that can be skipped.  The list could be
stored through a variety of methods (such as a hash table or binary search tree), although storing it as a simple binary array provides
for quicker access.  A binary array requires a number of elements equal to the modulus
product.  With a sufficient number of primes in the modulus 
product even this table can become large, in which case the modulus product can be split into 
different moduli (into multiple lookup tables), with each table able to reduce possibilities 
(unlike the a\_base\_skipahead table which must be just a single table).  Even though
Table~\ref{table-moduli} shows only a single a\_base\_elimination modulus for each exponent set,
for the purposes of code implementation, the moduli for exponent sets \{4,4,3\} and \{4,4,5\} permutation
cz\_minus\_ax were each split into three to provide three distinct a\_base\_elimination lookup tables.
And the moduli for \{5,5,3\} and \{5,5,4\} were each split into two to provide two distinct
a\_base\_elimination lookup tables.

Table~\ref{table1} shows the percentage of $c$ base and $f$ coefficients combinations eliminated for 
all $c$ bases and for $f$ coefficients ranging from 2 to 100,000.  
The elimination percentages could have been calculated based on all possible $f$ coefficients
(by ranging the coefficient residue from 0 to modulus-1), but then it would be impossible to exclude
coefficients that were not $z$-th power free.  The current table is more representative of actual
searching scenarios, and the ratios between the two calculation methods differ by only a few thousandths of a percent.

\section{Calculation of $f$ coefficient, $c$ base, and $a$ base Ranges}
\label{appendix-coeff}

This appendix details the logic to calculate the valid range of values for the loops in 
Algorithm~\ref{algorithm-Pegg} in \refsec\ref{subsection5dot-Pegg-algo}.
Let $V$ = the desired minimum Pegg Value.
For exponent sets \{$x,x,z$\}, the absolute minimum for the three bases is established 
by \refsec\ref{subsection1dot-smallest-bases} which states that the Pegg Value can be no higher than the 
smallest original base to the highest exponent.
For the exponent set \{3,4,5\}, Lemma~\ref{lemma-345} necessitates
the use of only the three coefficient possibilities listed below (and also provides that the
resultant equation will not be further multiplied by a number to 
re-associate $\min(A,B,C)$ to a base with a higher exponent).  For these \{3,4,5\} scenarios,
\refsec\ref{subsection1dot-smallest-bases} 
details how to determine the absolute minimum for the three bases.
We summarize this information here:
\begin{algorithmic}
\IF {exponent set = \{$x,x,z$\}, $z>x$}
  \STATE $a_{min1} \leftarrow 1$
  \STATE $b_{min1} \leftarrow 1$
  \STATE $c_{min1} \leftarrow V$
\ELSIF {exponent set = \{$x,x,z$\}, $z<x$}
  \STATE $a_{min1} \leftarrow V$
  \STATE $b_{min1} \leftarrow V$
  \STATE $c_{min1} \leftarrow 1$
\ELSIF {exponent set = \{3,4,5\}}
  \IF {coefficient is to be associated with the base to the 3rd power}
    \STATE [base to the 3rd power]$_{min1} \leftarrow V / 8$
    \STATE [base to the 4th power]$_{min1} \leftarrow V / 2$
    \STATE [base to the 5th power]$_{min1} \leftarrow V$
  \ELSIF {coefficient is to be associated with the base to the 4th power}
    \STATE [base to the 3rd power]$_{min1} \leftarrow V / 4$
    \STATE [base to the 4th power]$_{min1} \leftarrow V / 2$
    \STATE [base to the 5th power]$_{min1} \leftarrow V$
  \ELSIF {coefficient is to be associated with the base to the 5th power}
    \STATE [base to the 3rd power]$_{min1} \leftarrow V / 2$
    \STATE [base to the 4th power]$_{min1} \leftarrow V$
    \STATE [base to the 5th power]$_{min1} \leftarrow V$
  \ENDIF
\ENDIF
\end{algorithmic}
The base minimums may be further increased by considering the permutation type and combination of other
bases.

\subsection{Calculating $c$ base Range}
\label{subsectioncdot-cbase}

For a given $f$ coefficient,
the following logic is used to determine the valid minimum and maximum $c$ bases that
allow for the resultant equation to fit inside the desired range.

Let $N$ = the smallest equation multiplier (based on $f$) to convert the Original Pegg Equation to a resultant equation

Let $S_{min}$ = the desired minimum resultant equation size

Let $S_{max}$ = the desired maximum resultant equation size

\subsubsection{Permutation Type ax\_minus\_cz, $c$ base Minimum}
\label{subsubsectioncdot-ax-cz-c-min}

\begin{align*}
          a^x - f c^z &= b^y
\\
          f c^z + b^y &= a^x
\\
      N (f c^z + b^y) &= N a^x = \text{resultant equation size}
\\  
N (f c_{min2}^z + b^y) &\geq S_{min}
\intertext{
This information does not increase the minimum possible value for $c$
because $b$ can always be large enough
to make $N (f c_{min2}^z + b^y) \geq S_{min}$.
The minimum $c$ base is $c_{min1}$ as set initially.
}
  c_{min} &= c_{min1}
\end{align*}

\subsubsection{Permutation Type ax\_minus\_cz, $c$ base Maximum}
\label{subsubsectioncdot-ax-cz-c-max}

\begin{align*}
                   a^x- f c^z &= b^y
\\
                  f c^z + b^y &= a^x
\\
              N (f c^z + b^y) &= N a^x = \text{resultant equation size}
\\
N (f c_{max1}^z + b_{min1}^y) &\leq S_{max}
\\
   f c_{max1}^z + b_{min1}^y  &\leq S_{max}/N
\\
                 f c_{max}^z  &\leq S_{max}/N - b_{min1}^y
\\
                   c_{max}^z  &\leq (S_{max}/N - b_{min1}^y) / f
\\
                      c_{max} &= \bigg\lfloor \sqrt[z]{(S_{max}/N - b_{min1}^y) / f} \bigg\rfloor
\end{align*}

\subsubsection{Permutation Type cz\_minus\_ax, $c$ base Minimum}
\label{subsubsectioncdot-cz-ax-c-min}

For a given $f$, $c$ must be large enough so that the resulting equation will be $\geq S_{min}$:
\begin{align*}
       f c^z - a^x &= b^y
\\
             f c^z &= a^x + b^y
\\
           N f c^z &= N (a^x + b^y) = \text{resultant equation size}
\\
    N f c_{min2}^z &\geq S_{min}
\\
        c_{min2}^z &\geq S_{min} / (N f)
\\
          c_{min2} &= \Big\lceil \sqrt[z]{S_{min} / (N f)} \Big\rceil
\\
\intertext{
and $c$ must be large enough so that $f c^z \geq$ the sum of the minimum $a^x$ and $b^y$:
}
 f c^z - a^x &= b^y
\\
       f c^z &= a^x + b^y
\\
       f c^z &\geq a_{min1}^x + b_{min1}^y
\\
  c_{min3}^z &\geq (a_{min1}^x + b_{min1}^y) / f
\\
    c_{min3} &= \bigg\lceil \sqrt[z]{(a_{min1}^x + b_{min1}^y) / f} \bigg\rceil
\\
\intertext{
The minimum $c$ base is the maximum of the above two minimums and the $c_{min1}$ as set initially.
}
  c_{min} &= \max(c_{min1},c_{min2},c_{min3})
\end{align*}

\subsubsection{Permutation Type cz\_minus\_ax, $c$ base Maximum}
\label{subsubsectioncdot-cz-ax-c-max}

\begin{align*}
    f c^z - a^x &= b^y
\\
          f c^z &= a^x + b^y
\\
        N f c^z &= N (a^x + b^y) = \text{resultant equation size}
\\
   N f c_{max}^z &\leq S_{max}
\\
       c_{max}^z &\leq S_{max} / (N f)
\\
         c_{max} &= \Big\lfloor \sqrt[z]{S_{max} / (N f)} \Big\rfloor
\end{align*}

\subsection{Calculating $a$ base Range}
\label{subsectioncdot-abase}

For a given $f$ coefficient and $c$ base,
the following logic is used to determine the valid minimum and maximum $a$ bases that
allow for the resultant equation to fit inside the desired range.

Let $N$ = the smallest equation multiplier (based on $f$) to convert the Original Pegg Equation to a resultant equation

Let $S_{min}$ = the desired minimum resultant equation size

Let $S_{max}$ = the desired maximum resultant equation size

\subsubsection{Permutation Type ax\_minus\_cz, $a$ base Minimum}
\label{subsubsectioncdot-ax-cz-a-min}

For a given $f$ and $c$, $a$ must be large enough so that the resulting equation will be $\geq S_{min}$:
\begin{align*}
       a^x - f c^z &= b^y
\\
               a^x &= f c^z + b^y
\\
             N a^x &= N (f c^z + b^y) = \text{resultant equation size}
\\
      N a_{min2}^x &\geq S_{min}
\\
        a_{min2}^x &\geq S_{min} / N
\\
          a_{min2} &= \Big\lceil \sqrt[x]{S_{min} / N} \Big\rceil                          
\intertext{
and $a$ must be large enough so that $a^x \geq$ the sum of $f c^z$ and the minimum $b^y$:
}
       a^x - f c^z &= b^y
\\
               a^x &= f c^z + b^y
\\
        a_{min3}^x &\geq f c^z + b_{min1}^y
\\
          a_{min3} &= \bigg\lceil \sqrt[x]{f c^z + b_{min1}^y} \bigg\rceil
\intertext{
The minimum $a$ base is the maximum of the above two minimums and the $a_{min1}$ as set initially.
}
  a_{min} &= \max(a_{min1},a_{min2},a_{min3})
\end{align*}

\subsubsection{Permutation Type ax\_minus\_cz, $a$ base Maximum}
\label{subsubsectioncdot-ax-cz-a-max}

\begin{align*}
                a^x - f c^z &= b^y
\\
                        a^x &= f c^z + b^y
\\
                      N a^x &= N (f c^z + b^y) = \text{resultant equation size}
\\
                 N a_{max}^x &\leq S_{max}
\\
                   a_{max}^x &\leq S_{max} / N
\\
                   a_{max} &= \Big\lfloor \sqrt[x]{S_{max} / N} \Big\rfloor
\end{align*}

\subsubsection{Permutation Type cz\_minus\_ax, $a$ base Minimum}
\label{subsubsectioncdot-cz-ax-a-min}

For a given $f$ and $c$, $a$ must be large enough so that the resulting equation will be $\geq S_{min}$:
\begin{align*}
    f c^z - a^x &= b^y
\\
      a^x + b^y &= f c^z 
\\
  N (a^x + b^y) &= N (f c^z) = \text{resultant equation size}
\\
  N (a^x + b^y) &\geq S_{min}
\intertext{
This information does not increase the minimum possible value for $a$
because $b$ can always be large enough
to make $N (a^x + b^y) \geq S_{min}$.
With this permutation type $a^x + b^y = f c^z.$
If $x \neq y$, then
this information does not increase the minimum possible value for $a$
because $b$ can always be large enough
to make $a^x + b^y = f c^z$.  In this case:
}
  a_{min} &= a_{min1}
\intertext{
But when $x=y$ one of $a^x$ and $b^y$
must be $\geq f c^z / 2$, and the other must be $\leq f c^z / 2$.  
As established in \refsec\ref{subsection5dot-Pegg-algo}, for searching efficiency we 
label the larger value $a^x$ and the smaller value $b^y$.
When $x=y$, this increases the minimum $a$ base to:
}
          a_{min2}^x &\geq f c^z / 2
\\
            a_{min2} &= \Big\lceil \sqrt[x]{f c^z / 2} \Big\rceil
\\
             a_{min} &= \max(a_{min1},a_{min2})
\end{align*}

\subsubsection{Permutation Type cz\_minus\_ax, $a$ base Maximum}
\label{subsubsectioncdot-cz-ax-a-max}

\begin{align*}
      f c^z - a^x &= b^y
\\
              a^x &= f c^z - b^y
\\
        a_{max}^x &\leq f c^z - b_{min1}^y
\\
          a_{max} &= \bigg\lfloor \sqrt[x]{f c^z - b_{min1}^y} \bigg\rfloor
\end{align*}

\subsection{Calculating Valid $f$ coefficients}
\label{subsectioncdot-coeff}

For the exponent set \{3,4,5\}, Lemma~\ref{lemma-345} details the valid coefficients to be checked.
This section establishes the valid coefficients for exponent sets \{$x,x,z$\}.

Per \refsec\ref{subsection1dot-smallest-bases},
the Pegg Value of a resultant equation can be no greater 
than the smallest original base to the highest exponent.
We can use this information to determine the largest original coefficient
that could result in an equation having a minimum desired Pegg Value.
Given exponent set \{$x,x,z$\},

Let $V$ = the desired minimum Pegg Value

Let $S_{max}$ = the desired maximum resultant equation size

Let $R$ = the resultant coefficient as part of a resultant base to the highest exponent

Let $H$ = max($x,z$)

\begin{align*}
  (R_{max} V)^{H} &\leq S_{max}
\\
   R_{max} V &\leq \sqrt[H]{S_{max}}
\\
     R_{max} &\leq \sqrt[H]{S_{max}} / V
\\
     R_{max} &= \Big\lfloor \sqrt[H]{S_{max}} / V \Big\rfloor
\end{align*}

If $R > R_{max}$, it would limit the original base associated with the highest exponent 
to a value smaller than the desired Pegg Value, thereby precluding the resultant equation 
from having a Pegg Value $\geq V$.
But our search algorithm does not directly concern itself with resultant coefficients.  It loops 
through a range of original coefficients.  We must answer the question
``What is the highest original coefficient that would generate an $R \leq R_{max}$?''

Consider an Original Pegg Equation \{$x,x,z$\} with a single coefficient $f$ associated with the base to the $z$-th power.
With \{$x,x,z$\}, $z < x$, the two bases with the highest exponent will have a resultant coefficient 
that is based solely on the multiplier $N$. 
With \{$x,x,z$\}, $z > x$, the sole base with the highest exponent will have a resultant coefficient 
that is based on the product of $f$ and the multiplier $N$.
\\
Define the function cvt($x,z,v_p(f)$) for 
$x \geq 3;z \geq 3;$ 
$0 < v_p(f) < z$

if $z < x,$ cvt = $\smallest_{q}(v_p(f)) / \max(x,z)$

if $z > x,$ cvt = $(\smallest_{q}(v_p(f))+v_p(f)) / \max(x,z)$
\\
using the function:

$\smallest_{q}(v_p(f)) =$ smallest $q$ where $q \equiv 0 \pmod {x}, q \equiv -v_{p}(f) \pmod{z}$
\\
For a $p \mid f,$ cvt($x,z,v_p(f)$) returns the power of $p$ as represented in 
the resultant coefficient of the base(s) to the highest power.

\begin{table}[H]
\caption{
Power of $p$ as represented in $R$
(the resultant coefficient of the base(s) to the highest power).
Output of cvt function for all exponent sets \{$x,x,z$\} with each exponent $\leq 5.$
A subset of the information in Table~\ref{table2}.
\label{table-coeff2}}
\centering
\begin{tabular}{ c c c r r r r }
\hline\noalign{\smallskip}
         & Original     &          & & & &\\
         & Coefficient  &          & & & &\\
         & associated   &          & & & &\\
Exponent & with base to & Highest  & \multicolumn{4}{ c }{$v_{p}$(coeff)}\\
set      & this exponent& exponent & 1 & 2 & 3 & 4\\
\noalign{\smallskip}
\hline
\{4,4,3\} & 3 & 4 &  2 &  1 & \cellcolor[gray]{0.5} & \cellcolor[gray]{0.5} \\
\hline
\{5,5,3\} & 3 & 5 &  1 &  2 & \cellcolor[gray]{0.5} & \cellcolor[gray]{0.5} \\
\hline
\{3,3,4\} & 4 & 4 &  1 &  2 &  3 & \cellcolor[gray]{0.5} \\
\hline
\{5,5,4\} & 4 & 5 &  3 &  2 &  1 & \cellcolor[gray]{0.5} \\
\hline
\{3,3,5\} & 5 & 5 &  2 &  1 &  3 &  2 \\
\hline
\{4,4,5\} & 5 & 5 &  1 &  2 &  3 &  4 \\
\hline
\end{tabular}
\end{table}

With the cvt function to show how a prime dividing an original coefficient is
represented in the resultant coefficient of the base(s) to the highest exponent,
we can answer the question
``What is the highest original coefficient that would generate an $R \leq R_{max}$?''
The highest possible original coefficient must be 
$\leq {R_{max}}^T$ where $T$ = maximum ratio of $v_p(f)$ to cvt($x,z,v_p(f)$).

Consider the exponent set \{5,5,4\}.  For $p \mid f$ with $v_p(f) = 3$, the entire equation must be multiplied
by $p^{5}$ (as shown in Table~\ref{table1}).  This is represented as part of the resultant coefficient of
the bases to the highest power as $p^{5 / z}$ and the cvt function shows this
with the return value of 1 (representing $\log_p(p^{5 / 5})$.
So with $v_p(f)=3$, an original coefficient of $p^{3}$ generates an $R$ of $p$.  
This 3:1 power ratio is the highest ratio for the \{5,5,4\} exponent set.
To produce an $R \leq R_{max},$ the original coefficient of a \{5,5,4\} equation must be $\leq {R_{max}}^3$.

For \{$x,x,z$\} with all exponents $\leq 5$, \{5,5,4\} has a highest ratio of 3, \{4,4,3\} and \{3,3,5\} each have
a highest ratio of 2, and the remaining exponent sets \{5,5,3\}, \{3,3,4\}, \{4,4,5\} each have a highest ratio
of 1.

At this point, the list of possible original coefficients is all integers $\geq 2$ and $\leq {R_{max}}^T$. 
The strategy is to create a ``valid original coefficient'' boolean array of these integers, initially 
marking all the entries as valid.  Then we will perform several steps, each time possibly ruling out
candidates.
\\

Step 1: Per \refsec\ref{subsection1dot-powerfree}, mark as invalid any coefficients that are not $z$-th power free.
\\

Step 2: Mark as invalid any coefficients that have a smallest multiplier $N$ that is sufficiently large that it
precludes an original equation that has both:
\begin{compactitem}
\item{
a base to the highest exponent that is $\geq$ the desired Pegg Value
}
\item{
that converts to a resultant equation $\leq$ the maximum equation size
}
\end{compactitem}
  
For example, exponent set \{3,3,5\} requires a coefficient 
with $v_{p}$($f$)$=2$ to have the entire equation multiplied by $p^3$, but a coefficient with 
$v_{p}$($f$)$=1$, 3, or 4 requires a larger multiplier.  
Consider the search for Pegg Values $> 63742$ in \{3,3,5\} equations $\leq 2^{88}$
The highest possible resultant coefficient is 3, as
\begin{align*}
     R_{max} &= \Big\lfloor \sqrt[H]{S_{max}} / V \Big\rfloor
\\
     R_{max} &= \Big\lfloor \sqrt[5]{2^{88}} / 63743 \Big\rfloor
\\
     R_{max} &= 3
\intertext{
Let $T$ = 2 = the maximum \{3,3,5\} ratio of $v_p(f)$ to cvt($x,z,v_p(f)$).
The maximum original coefficient is given by:
}
     O_{max} &\leq {R_{max}}^T
\\
     O_{max} &\leq 3^2
\\
     O_{max} &\leq 9
\end{align*}

Examining each possible coefficient $\geq 2$ and $\leq 9$, only the coefficients 4 and 9 can product a resultant
equation $\leq 2^{88}$.  Step 2 marks as invalid the coefficients 2, 3, 5, 6, 7, 8.
\\

Step 3: for each coefficient still marked as valid, for the given permutation type,
calculate the minimum and maximum possible $c$ bases
that could result in an equation with a desired Pegg Value and within the resultant equation range.
If a coefficient produces a minimum $c$ base $>$ its maximum $c$ base, the coefficient is marked as invalid.
For example,
consider the search for Pegg Values $> 63742$ in \{5,5,3\} equations $\leq 2^{100}$ of permutation type
cz\_minus\_ax.
The highest possible resultant coefficient is 16, as
\begin{align*}
     R_{max} &= \Big\lfloor \sqrt[H]{S_{max}} / V \Big\rfloor
\\
     R_{max} &= \Big\lfloor \sqrt[5]{2^{100}} / 63743 \Big\rfloor
\\
     R_{max} &= 16
\intertext{
The coefficient 16 is marked as invalid in Step 1 because it is not $z$-th power free.
}
\intertext{
For the original coefficient 15, the minimum $c$ base is given by \refsec\ref{subsubsectioncdot-cz-ax-c-min}:
}
c_{min} = c_{min3} &= \bigg\lceil \sqrt[z]{(a_{min1}^x + b_{min1}^y) / f} \bigg\rceil
\\
c_{min} = c_{min3} &= \Big\lceil \sqrt[3]{(63743^5 + 63743^5) / 15} \Big\rceil
\\
c_{min} = c_{min3} &= 51963742
\intertext{
For the original coefficient 15, the maximum $c$ base is given by \refsec\ref{subsubsectioncdot-cz-ax-c-max}:
}
c_{max} &= \Big\lfloor \sqrt[z]{S_{max} / (N f)} \Big\rfloor
\\
c_{max} &= \Big\lfloor \sqrt[3]{2^{100} / ((3^5 * 5^5) * (3*5))} \Big\rfloor
\\
c_{max} &= 48100619
\end{align*}
So coefficient 15 is excluded as it has no valid $c$ bases.  All the other coefficients that passed the tests 
in step 1 and step 2 are valid, as each coefficient has a minimum $c$ base $\leq$ its maximum $c$ base.
\\
\\
After step 3, all the coefficients still marked as valid are to be checked.

\end{document}